\documentclass[reqno,draft,12pt]{amsart}
\usepackage{amsmath}
\usepackage{amscd}
\usepackage{amssymb}
\usepackage{enumerate}
\usepackage{amsfonts}
\usepackage{graphicx}
\usepackage{latexsym}
\usepackage[all]{xy}
     
\numberwithin{equation}{section}

\newcommand{\op}[1]{\operatorname{#1}}

\newcommand{\Gwh} {\widehat{G}}
\newcommand{\Fwh} {\widehat{F}}
\newcommand{\Fiwh} {{\widehat{F}}^{-1}}
\newcommand{\cKwh} {\widehat{\cK}}
\newcommand{\rar}[1]{\stackrel{#1}{\longrightarrow}}
\newcommand{\lar}[1]{\stackrel{#1}{\longleftarrow}}
\newcommand{\myiota}{\mbox{{\Large ${\iota}$}}}

\newcommand{\al}{\alpha}
\newcommand{\be}{\beta}
\newcommand{\Betawh} {\widehat{\beta}}

\newcommand{\om}{\omega}

\newcommand{\bC}{{\mathbb C}}
\newcommand{\bR}{{\mathbb R}}

\newcommand{\cP}{{\mathcal P}}
\newcommand{\cQ}{{\mathcal Q}}

\newcommand{\cJ}{{\mathcal J}}

\newcommand{\cK}{\underline{\mathcal J}}

\newcommand{\cY}{{\mathcal Y}}
\newcommand{\cX}{Z}
\newcommand{\cXwh}{\widehat{Z}}

\newcommand{\cL}{{\mathcal L}}
\newcommand{\cM}{{\mathcal M}}
\newcommand{\cS}{{\mathcal S}}

\newcommand{\cF}{{\mathcal F}}

\newcommand{\Dwh} {\widehat{D}}

\newcommand{\ph} {\hat{p}}
\newcommand{\qh} {\hat{q}}

\newcommand{\pih} {\hat{\pi}}
\newcommand{\alh} {\hat{\alpha}}
\newcommand{\gh} {\hat{g}}
\newcommand{\bh} {\hat{b}}
\newcommand{\hh} {\hat{h}}
\newcommand{\sh} {\hat{s}}
\newcommand{\yh} {\hat{y}}
\newcommand{\Bwh} {\widehat{B}}

\newcommand{\Jwh} {\widehat{J}}
\newcommand{\Lwh} {\widehat{L}}
\newcommand{\jh} {\hat{j}}

\newcommand{\Xwh} {\widehat{X}}
\newcommand{\Ywh} {\widehat{Y}}
\newcommand{\Zwh} {\widehat{Z}}
\newcommand{\cJwh} {\widehat{\cJ}}

\newcommand{\cJpwh} {\widehat{\cJ}'}
\newcommand{\cKpwh} {\widehat{\cK}'}

\newcommand{\pdxii} {\partial/{\partial \xi_i}}
\newcommand{\pdxij} {\partial/{\partial \xi_j}}
\newcommand{\pdxi} {\partial/{\partial x_i}}
\newcommand{\pdyi} {\partial/{\partial y_i}}
\newcommand{\smo}{i}

\newcommand{\Aut}{\operatorname{Aut}}

\newcommand{\Wedge}{\bigwedge}

\renewcommand{\dot}{^{\bullet}}

\newcommand{\sbr}{\smallbreak}

\newcommand{\Ann}{\operatorname{Ann}}

\newcommand{\cou}[2]{\bigl[#1,#2\bigr]}

\newcommand{\bra}{\langle}
\newcommand{\ket}{\rangle}
\newcommand{\pairing}{\bra\bullet,\bullet\ket}
\newcommand{\pair}[2]{\bra #1, #2\ket}
\newcommand{\eval}[2]{\bra #1\,\big\vert\,#2\ket}

\newcommand{\matr}[4]{\begin{pmatrix}  #1 \qquad &  #2  \\ 
    #3 \qquad & #4  \end{pmatrix}}

\newtheorem{thm}{Theorem}[section]
\newtheorem{cor}[thm]{Corollary}
\newtheorem{lem}[thm]{Lemma}
\newtheorem{prop}[thm]{Proposition}

\theoremstyle{remark}
\newtheorem{rem}[thm]{Remark}

\newtheorem{example}[thm]{Example}
\newtheorem{examples}[thm]{Examples}
\newtheorem{defin}[thm]{Definition}
\newtheorem{conj}[thm]{Conjecture}

\newlength{\dummy}

\title[Mirror Symmetry and Generalized Complex Manifolds]{Mirror
  Symmetry and Generalized Complex Manifolds}  

\author[O.~Ben-Bassat]{Oren Ben-Bassat\address{Oren Ben-Bassat:
Department of Mathematics, University of Pennsylvania,
Philadelphia, PA 19104-6395, e-mail: orenb@math.upenn.edu}}

\date{}

\begin{document}

\begin{abstract} 
In this paper we develop a relative version of T-duality in
generalized complex geometry which we propose as a manifestation of
mirror symmetry.  Let $M$ be an $n-$dimensional smooth real manifold,
$V$ a rank $n$ real vector bundle on $M$, and $\nabla$ a flat
connection on $V$.  We define the notion of a $\nabla-$semi-flat
generalized complex structure on the total space of $V$.  We show that
there is an explicit bijective correspondence between
$\nabla-$semi-flat generalized complex structures on the total space
of $V$ and $\nabla^{\vee}-$semi-flat generalized complex structures on
the total space of $V^{\vee}$.  Similarly we define semi-flat
generalized complex structures on real $n-$torus bundles with section
over an $n$-dimensional base and establish a bijective
correspondence between semi-flat generalized complex structures on
pairs of dual torus bundles.  Along the way, we give methods of
constructing generalized complex structures on the total spaces of
vector bundles and torus bundles with sections.  We also show that
semi-flat generalized complex structures give rise to a pair of
transverse Dirac structures on the base manifold.  We give
interpretations of these results in terms of relationships between the
cohomology of torus bundles and their duals.  We also study the ways
in which our results generalize some well established aspects of mirror
symmetry as well as some recent proposals relating 
generalized complex geometry to string theory.
\end{abstract}

\maketitle

\tableofcontents

\section{Introduction}\label{s:intro}

Mirror symmetry is often thought of as relating the very different
worlds of complex geometry and symplectic geometry.  It was recently
shown by Hitchin \cite{H2} that symplectic and complex structures on a
manifold have a simple common generalization called a generalized
complex structure.  This is a complexified version of Dirac geometry
\cite{Cou} along with an extra non-degeneracy condition.  It is
expected that Mirror Symmetry should give rise to an involution on
sectors of the moduli space of all Generalized Complex Manifolds of a
fixed dimension.  One of the most concrete descriptions of the mirror
correspondence is the Strominger-Yau-Zaslow picture \cite{SYZ} in
which Mirror Symmetry is interpreted as a relative T-duality along the
fibers of a special Lagrangian torus fibration. This is sometimes
referred to as ``T-Duality in half the directions''.  In our previous
work \cite{BBB}, we investigated the linear algebraic aspects of
T-duality for generalized complex structures.  See also \cite{WeiTan}
for the analogous story in Dirac geometry.  In this paper, we go one
step further and construct an explicit mirror involution on certain
moduli of generalized complex manifolds.  Similarly to the case of
Calabi-Yau manifolds the definition of our mirror involution depends
on additional data.  In our set up, we will consider generalized
complex manifolds equipped with a compatible torus fibration.  This
involution, when applied to such a manifold, gives another with the
same special properties, which we propose to identify as its mirror
partner.  In the special cases of a complex or symplectic structure on
a semi-flat Calabi-Yau manifold our construction reproduces the
standard T-duality of \cite{Pol,Leung,Pioli}. In addition we get 
new examples of mirror symmetric generalized complex manifolds,
e.g. the ones coming from B-field transforms of complex or symplectic
structures.

If $V$ is a real vector space then \cite{H2} a generalized complex
structure on $V$ is a complex subspace $E \subseteq (V \oplus V^{\vee})
\otimes \bC$ that satisfies $E\cap\overline{E}=(0)$ and is
maximally isotropic with respect to the canonical quadratic form on
$(V \oplus V^{\vee}) \otimes \bC$. Let
\[
f : V \oplus V^{\vee} \to W\oplus
W^{\vee}
\] 
be a linear isomorphism which is compatible with the
canonical quadratic forms.  Then $f$ induces a bijection between
generalized complex structures on $V$ and generalized complex
structures on $W$.  Transformations of this type can be viewed as
linear analogues of the T-duality transformations investigated in the
physics literature (see \cite{K,Van} and references therein).
Mathematically they were studied in \cite{WeiTan} for Dirac
structures and in \cite{BBB} for generalized complex structures.  
In this paper, the relevant case is where $V = A \oplus B$, $W =
A^{\vee} \oplus B$, and $f : A \oplus B \oplus A^{\vee} \oplus
B^{\vee} \to A^{\vee} \oplus B \oplus A \oplus B^{\vee}$ is the
obvious shuffle map.

A generalized complex structure on a manifold $X$ is a maximally isotropic 
sub-bundle of $(T_X \oplus T_X^{\vee})\otimes \bC$ that satisfies 
$E\cap\overline{E}=(0)$ and that $E$ is closed under the Courant Bracket.  
In this paper, we shall preform a relative version of this T-duality for 
pairs of manifolds that are fibered over the same base and where the 
two fibers over each point are ``dual'' to each other.  In other words we 
will find a way to apply the linear ideas above to the torus fibered approach.  On each fiber, this process will agree with the linear map described above.  

We relate the integrability of {\em semi-flat} (see definition 
\ref{def:semiflat})  
generalized almost complex  
structures on torus and vector bundles to data which lives only on the 
base manifold.  We show that a semi-flat generalized almost complex 
structure is integrable if and only its mirror structure is integrable.

Using a natural connection on a torus bundle $Z \to M$ with 
zero section $s$, we will construct 
semi-flat generalized complex structures $\cJ$ on $\cX$ from generalized 
almost complex structures $\cK$ on the vector bundle 
$s^*T_{\cX/M} \oplus T_M$.  The definition of semi-flat includes 
the condition that 
\[
\cK(s^*T_{\cX/M} \oplus s^*T_{\cX/M}^{\vee}) = T_M \oplus T_M^{\vee}.
\]

\noindent
Then we have the following two results:

\begin{thm} (\ref{thm:torintegrability})
A semi-flat generalized almost complex structure $\cJ$ on a torus bundle 
$\cX \to M$ with zero section $s$ is integrable if and only if 
\[
\bigl[\cK (\cS \oplus \cS^{\vee}),\cK (\cS \oplus \cS^{\vee}) \bigr] = 0
\]  
where $\cS$ is the sheaf of flat sections of $s^*T_{\cX/M}$.
\end{thm}

\begin{cor} (\ref{cor:tormirror})
A semi-flat generalized almost complex structure $\cJ$ on a torus bundle 
$\cX \to M$ is integrable if and only if its mirror structure $\cJwh$ on the 
dual torus bundle $\cXwh \to M$ is integrable.
\end{cor}

\

\noindent
These statements set the stage for understanding mirror symmetry and
the mirror transform of D-branes in generalized Calabi-Yau
geometry. Our results are a direct generalization of the setup employed by 
Arinkin and Polishchuk \cite{Pol} in ordinary mirror 
symmetry.  Explicit 
examples of this fact can be found in 
section \ref{s:examples}.

We relate this transformation of geometric structures to a purely
topological map on differential forms which descends to a map from the
de Rham cohomology of $Z$ to the de Rham cohomology of $\Zwh$.  In
particular, the map on differential forms exchanges the pure
spinors associated to the generalized complex structure on $Z$ with
the ones associated to the mirror generalized complex structure on
$\Zwh$.  This type of transformation was also discussed in \cite{Allessandro}.

Throughout the paper we comment on how our results relate to some of
the well established results and conjectures of mirror 
symmetry \cite{Pol, SYZ, Leung, Pioli, Ko} and also what they 
say in regards to the new developments in
generalized K\"{a}hler geometry \cite{Gua} and the relationships
between generalized complex geometry and string theory \cite{K, Van,
Gua} which have appeared recently.  As mentioned in \cite{K} we may
interpret these dualities as being a generalization of the duality
between the $A-$model and $B-$model in topological string theory.  In
the generalized K\"{a}hler case, they can be interpreted as dualities
of supersymmetric nonlinear sigma models \cite{Gates}.  To this end,
in section~\ref{s:Branes} we sketch a relationship between branes in
the sense of \cite{Gua, K} in a semi-flat generalized complex
structure and branes in its mirror structure.  For some simple examples 
of branes, we give the relationship directly.  We also show in 
section~\ref{s:vectmirror} that the 
Buscher rules \cite{Bus1, Bus2} for the transormation of metric and 
B-field hold between the mirror pairs of generalized K\"{a}hler manifolds 
that we consider.

It will be very interesting to extend the discussion in
section~\ref{s:Branes} to a full-fledged Fourier-Mukai transform on
generalized complex manifolds.  Unfortunately, the in-depth study of
branes in generalized complex geometry is obstructed by the
complexity of the the behavior of sub-manifolds with regards to a
generalized complex structure.  Several subtle issues of this nature
were analyzed in our previous paper \cite{BBB}. In particular we
investigated in detail the theory of sub and quotient generalized
complex structures, described a zoo of sub-manifolds of generalized
complex manifolds and studied the relations among those. We also gave
a classification of linear generalized complex structures and
constructed a category of linear generalized complex structures which
is well adapted to the question of quantization. In a future work we
plan to incorporate the structure of a torus bundle in this analysis
and construct a complete Fourier transform for branes.

For the benefit of the reader who may not be familiar with generalized 
complex geometry, we have included \S\ref{s:notation} which introduces the 
linear algebra and some basics on generalized complex 
manifolds.  More details on these basics may be found 
in \cite{BBB, Gua, H2, H3, Huy}.

\section{Notation, conventions, and basic definitions}\label{s:notation}

Overall, we will retain the notation and conventions from our previous 
paper \cite{BBB}, and so we only recall the most important facts for this 
paper as well as some changes.  The dual of a vector space 
$V$ will be denoted 
as $V^{\vee}$.  We will 
often use the annihilator of a subspace $W\subseteq V$, which we will denote

\[
\Ann(W)=\left\{f\in V^{\vee}\left| f\big\vert_W\equiv 0\right.\right\} \subseteq V^{\vee}.
\]

\noindent
We will need the pairing $\pairing$ on
$V\oplus V^{\vee}$, given by (following \cite{H3})
\[
\pair{v+f}{w+g}=-\frac{1}{2}\bigl(f(w)+g(v)\bigr) \text{  for all
} v,w\in V,\  f,g\in V^{\vee}.
\]

\noindent
Given $v\in V$ and $f\in V^{\vee}$, we will write either
$\eval{f}{v}$ or $\eval{v}{f}$ for $f(v)$. This pairing 
corresponds to the quadratic form $Q(v+f)=-f(v)$.

\sbr

We will tacitly identify elements $B\in\Wedge^2 V^{\vee}$ with 
linear maps $V \to V^{\vee}$.  When thought of in this way, 
we have that the map is skew-symmetric: $B = -B^{\vee}$. 

We will often consider linear maps of $V\oplus W \to V' \oplus W'$. 
Sometimes, these be written as  matrices
\[
T=\matr{T_1}{T_2}{T_3}{T_4},
\]
with the understanding that $T_1:V\to W'$, $T_2:W\to V'$,
$T_3:V\to W'$ and $T_4:W\to W'$ are linear maps.  All of these 
conventions will be extended to vector bundles 
and their sections in the obvious way.
\sbr

If $M$ is a manifold, we let
$C^\infty_M$ denote the sheaf of {\em
real-valued} $C^\infty$ functions on $M$.  We will use the 
same notation for a vector bundle and for its sheaf of sections.  
The tangent and cotangent bundles of $M$ will be denoted by 
$T_M$ and $T_M^{\vee}$.  For a vector bundle $V$ over a manifold $M$ and a 
smooth map $f:N \to M$, we denote the pullback bundle by 
$f^*V$.  A section of $f^*V$ which is a pullback of a 
section $e$ of $V$ will be 
denoted $f^*(e)$.  If f is an isomorphism onto its image or the projection 
map of a fiber bundle, the sections of this form give the sub-sheaf 
$f^{-1}V \subseteq f^* V$.  We will sometimes replace  
$\Wedge\dot T^{\vee}_{M}$ by $\Omega_{M}\dot$. 

Now we will give some basic facts on generalized complex
geometry that we will need in the paper.  For more information the reader may 
see \cite{BBB, H2, Gua}. 

\subsection{Generalized almost complex manifolds}\label{ss:gacm}

Let $M$ be a real manifold. A {\em generalized almost complex
structure} on a real vector bundle $V\to M$ has been defined 
\cite{Gua, H2,H3} in the following equivalent ways:

\begin{itemize}
\item A sub-bundle $E\subseteq V_{\bC}\oplus V_{\bC}^{\vee}$ which is
maximally isotropic with respect to the standard pairing
$\pairing$ and satisfies $E\cap\overline{E}=0$
\item An automorphism $\cJ$ of $V\oplus V^{\vee}$
which is orthogonal with respect to $\pairing$ and satisfies
$\cJ^2=-1$.
\end{itemize}

\begin{example}\label{d:complsympl} 
Let $V$ be a real vector bundle.
\begin{enumerate}[(a)]
\item Let $J$ be an almost complex structure on $V$. Then
\[
\cJ=\matr{J}{0}{0}{-J^{\vee}}
\]
\noindent
is a generalized almost complex structure on $V$. If $\cJ$ is a generalized
complex structure on $V$ that can be
written in this form, we say that $\cJ$ is of {\em complex type}.
\item Let $\om$ be an almost symplectic form on $V$ (i.e.,
a non-degenerate section $\om$ of $\Wedge^2 V^{\vee}$). Then
\[
\cJ=\matr{0}{-\om^{-1}}{\om}{0}
\]
\noindent
is a generalized complex structure on $V$. We say that such  a $\cJ$
is of  {\em symplectic type}.
\end{enumerate}
\end{example}

There is also a way of describing generalized almost complex 
structures on $V$ in terms 
of line sub-bundels of $\Wedge V^{\vee} \otimes \bC$ or spinors.  This 
interpretation is very convenient for some purposes. 

\begin{defin}\label{d:purespbun} \cite{H3, Gua} 
Let $\cJ$ be a generalized almost complex structure on 
a vector bundle $V$ over $M$.  Define   
the {\em canonical bundle} to be the complex line bundle 
$L \subseteq \Wedge^{\dot} V^{\vee}\otimes{\mathbb{C}}$ 
consist of the sections $\phi$ satisfying 
$\myiota_v \phi + \al \wedge \phi$ for all sections 
$v + \alpha$ of the $+\smo$ eigenbundle $E$ corresponding to the 
generalized almost complex structure on $V$.  Sections of $L$ will be called 
{\em representative spinors}.
\end{defin}

For the case of an almost symplectic manifold with two form $\om$, this 
line bundle is generated by $\exp(-i \om)$.  For an almost complex manifold, 
one gets the usual canonical bundle.  Spinor bundles can also be 
understood intrinsically in terms of the 
sheaves of modules over appropriate sheaf of Clifford algebras.  
The sections will satisfy certain restrictions over each fiber.  
They are known as pure spinors \cite{Che, Gua, H3}.  We have listed some 
of their features and examined their restriction to submanifolds in \cite{BBB}.

\begin{defin}\cite{H3}
In the special case that $V = T_M$ has a 
generalized almost complex structure, we call $M$ a {\em generalized
almost complex manifold}.
\end{defin}
In this case the spinor sections are 
differential forms.  Such a manifold is always even dimensional
as a real manifold.  This can be shown by constructing two
almost complex structures on $M$ out of the generalized almost complex
structure \cite{H3}.  This also follows from the classification of
generalized complex structures on a vector space which was done 
in our previous paper \cite{BBB}.  For the case of manifolds, a local 
structure theorem for generalized complex maniolds 
has been proven by Gualtieri \cite{Gua}.

Consider a real vector bundle $V$ and an automorphism 
$\cJ$ of $V\oplus V^{\vee}$, written in matrix form as

\[
\cJ=\matr{\cJ_1}{\cJ_2}{\cJ_3}{\cJ_4}.
\]
Let us record, for future use, the restrictions on the $\cJ_i$ coming from
the conditions that $\cJ$ preserves the pairing $\pairing$
and satisfies $\cJ^2=-1$.  They are:
\begin{eqnarray}
\label{e:1} &&\cJ_1^2+\cJ_2\cJ_3=-1; \\
\label{e:2} &&\cJ_1\cJ_2+\cJ_2\cJ_4=0; \\
\label{e:3} &&\cJ_3\cJ_1+\cJ_4\cJ_3=0; \\
\label{e:4} &&\cJ_4^2+\cJ_3\cJ_2=-1; \\
\label{e:5} &&\cJ_4=-\cJ_1^{\vee}; \\
\label{e:6} &&\cJ_2^{\vee}=-\cJ_2; \\
\label{e:7} &&\cJ_3^{\vee}=-\cJ_3.
\end{eqnarray}

\subsection{$B$- and $\be$-field transforms} \cite{H1,H2,H3,Gua} \label{ss:B-field}
Consider a real vector bundle $V$ and a global 
section $B$ of $\Wedge^2 V^{\vee}$.  
Consider 
the transformation of $V \oplus V^{\vee}$
\[
\exp(B):=\matr{1}{0}{B}{1}.
\]
It is easy to see that $\exp(B)$ is an orthogonal
automorphism of $V\oplus V^{\vee}$. Thus $\exp(B)\cdot E$ is a 
generalized almost complex structure on $V$ for any generalized almost 
complex structure $E \subseteq (V \oplus V^{\vee})\otimes{\bC}$ on $V$.  
We will call $\exp(B)\cdot E$ the {\em
$B$-field transform of $E$ defined by $B$}.  We should note here that 
these type of transformations are sometimes called gauge-transformations and 
were introduced with that name into real Dirac geometry \cite{Cou} 
in \cite{SW}.  For an overview of 
these transformations in the Dirac geometry context, see \cite{BR}.  
Similarly, if
$\be\in\Wedge^2 V$, then 
\[
\exp(\beta):=\matr{1}{\be}{0}{1}
\]
then $\exp(\beta)\cdot E$ will be called the {\em $\be$-field
transform of $E$ defined by $\be$}.  One can also write these transfomations 
in terms of the orthogonal automorphisms
$\cJ$ of $V\oplus V^{\vee}$.  In this case, the actions of $B$ and $\be$ are
given by $\cJ \mapsto \exp(B) \cJ \exp(-B)$ and $\cJ \mapsto \exp(\beta)\cJ
\exp(-\beta)$, respectively.  We can also describe $B$-field transforms
in terms of local spinor representatives: if a generalized almost complex 
structure on a real vector bundle 
$V$ is defined by a pure spinor $\phi\in\Wedge\dot V_{\bC}^{\vee}$,
and $B\in\Wedge^2 V^{\vee}$ then the $B$-field transform of this
structure corresponds to the pure spinor 
$\exp(-B)\wedge\phi$ \cite{H2,H3}.  The $\be$-field transform 
corresponds to the pure spinor 
$\myiota_{\exp(\be)}\phi$ \cite{Gua,H3}.

\subsection{Integrability}\label{ss:integrability} 
Let $M$ be a generalized almost complex manifold.  The {\em Courant
bracket} (\cite{Cou}, p.~645) is defined on sections of 
$(T_M \oplus T^{\vee}_M)\otimes{\bC}$ by 

\[
\bigl[{X + \xi},{Y + \eta}\bigr]=[X,Y]+\cL_X\eta-\cL_Y\xi+\frac{1}{2}\cdot
d(\myiota_Y\xi-\myiota_X\eta).
\]
or equivalently
\[
\bigl[{X + \xi},{Y + \eta}\bigr]=[X,Y]+\myiota_X d\eta +\frac{1}{2}d\myiota_X\eta -\myiota_Y d\xi-\frac{1}{2}\myiota_Y\xi.
\]

\begin{defin}[cf. \cite{Gua}, \cite{H2}, \cite{H3}]\label{d:gcm}
Let $M$ be a real manifold equipped with a generalized almost complex
 structure defined by $E\subseteq (T_M \oplus
 T^{\vee}_M)\otimes{\bC}$. We say that $E$ is {\em integrable} if the
 sheaf of sections of $E$ is closed under the Courant bracket. If that
 is the case, we also say that $E$ is a {\em generalized complex
 structure} on $M$, and that $M$ is a {\em generalized complex
 manifold}.
\end{defin}

\begin{rem}[cf. \cite{Cou},\cite{H3}]
As we noted in \cite{BBB}, the integrability condition for 
a generalized almost complex structure $\cJ$ is equivalent 
\cite{BBB} to the vanishing of the
{\em Courant-Nijenhuis} tensor
\[
N_{\cJ}(X,Y)=\cou{\cJ X}{\cJ Y}-\cJ\cou{\cJ X}{Y}-\cJ\cou{X}{\cJ
Y}-\cou{X}{Y}
\]
where $X$, $Y$ are sections of $T_M\oplus T_M^{\vee}$.
\end{rem}

Integrability can also be expressed in terms of spinors \cite{H3,Gua, Xu}.  
If $L \subseteq \Wedge^{\dot} T_M^{\vee} \otimes{\mathbb{C}}$ is the \
line bundle of spinors 
associated to 
a generalized almost 
complex structure $\cJ$ on a manifold $M$ then $\cJ$ is integrable if and 
only if all sections $\phi$ of $L$ satisfy 
\[
d\phi = \myiota_v \phi + \alpha \wedge \phi
\]
for some section $v + \alpha$ of $(T_M \oplus T_M^{\vee})\otimes{\mathbb{C}}$.

\
\begin{examples}[cf. \cite{H3, Gua}]
In the case that a generalized almost complex structure comes from an almost 
complex structure, it will be integrable if and only if the almost complex 
structure is integrable, giving a complex structure to the manifold.  
In the case that a generalized almost 
complex structure comes from a nondegenerate differential $2$-form 
(almost symplectic structure), it will 
be integrable if and only if the form is closed, i.e. gives a 
sympletic structure to the manifold.
\end{examples}

General $B-$field and $\beta-$field transformations 
need not preserve integrability.  However \cite{H2}, a {\em
closed} $2$-form $B$ acts on generalized complex structures on $M$ in
the same way as described in \S\ref{ss:B-field}.  In fact, 
a $B$-transform by
a $2$-form on $M$ is an automorphism of the Courant bracket if and
only if the $2$-form is closed \cite{H3}. Note furthur that a
$B$-field transform of a particular generalized complex structure can
be integrable even if the $2$-form is not closed. In fact, for any
specific generalized complex manifold $(M,\cJ)$, one can write
down explicitly the conditions that need to be satisfied by a two form, $B$
or a bi-vector field $\beta$ in order for the $B$-field or
$\beta$-field transform of $(M, \cJ)$ to be integrable.  We will study
examples of this phenomenon in Section \ref{s:examples}.

\subsection{Generalized almost K\"{a}hler manifolds} \label{ss:kahlerdef}
We will need the notion \cite{H3, Gua} of a generalized almost K\"{a}hler
structure.

\begin{defin} \cite{Gua} \label{def:genkahler}
A generalized almost K\"{a}hler structure on a manifold $M$ is specified by 
one of the equivalent sets of data.

1)  A pair $(\cJ,\cJ')$ of commuting generalized almost complex structures 
whose product, $G = -\cJ \cJ'$ is positive definite with respect to the
standard quadratic form $\pairing$ on $T_M \oplus T_M^{\vee}$.

2)  A quadruple $(g,b,J_{+},J_{-})$ consisting of a Riemanian metric $g$, 
two-form $b$, and two almost complex structures $J_{+}$ and $J_{-}$ such  
that the isomorphisms $\om_+ = g J_{+}: T_M \to T_M^{\vee}$ 
and $\om_- = g J_{-}: T_M \to T_M^{\vee}$ are anti-symmetric and hence 
correspond to non-degenerate two-forms. 

The two sets of data are related explicitly as follows.  The $(+1)$ 
eigenbundle of 
$G$ is the graph of $g + b : T_M \to T_M^{\vee}$.  Denote this vector bundle
by $C_{+}$, and the $(-1)$ eigenbundle (which is the graph of $b-g$) 
by $C_{-}$.  Then

\[
J_{\pm} = \pi_{T_M} \circ \cJ \circ {(\pi_{T_M} \big\vert_{C_{\pm}})}^{-1}.
\]

\noindent
Conversely, given $(g,b,J_{+},J_{-})$, one defines

\[
\cJ=\frac{1}{2}\left(\begin{matrix}1&0\\b&1\end{matrix}\right)
\left(\begin{matrix}J_+ + J_- & -(\omega_+^{-1} - \omega_-^{-1}) \\
\omega_+ - \omega_-&-(J^{\vee}_+ +  J^{\vee}_-)\end{matrix}\right)
\left(\begin{matrix}1&0\\-b&1\end{matrix}\right)
\]

\noindent
and

\[
\cJ'=\frac{1}{2}\left(\begin{matrix}1&0\\b&1\end{matrix}\right)
\left(\begin{matrix}J_+ - J_- & -(\omega_+^{-1} + \omega_-^{-1}) \\
\omega_+ + \omega_-&-(J^{\vee}_+ -  J^{\vee}_-)\end{matrix}\right)
\left(\begin{matrix}1&0\\-b&1\end{matrix}\right)
\]

\end{defin}

\

\noindent
Using this same notation we have that 

\begin{equation}\label{e:G}
G=\left(
\begin{array}{cc}
-g^{-1}b & g^{-1}  \\
g - bg^{-1}b & bg^{-1} 

\end{array}
\right).
\end{equation}

\begin{examples} \cite{Gua}  \label{exs:gualt}
Notice that this definition naturally generalizes the linear 
algebraic data of an 
K\"{a}hler manifold.  We will refer to this as the
ordinary K\"{a}hler case.  There is an important family of examples which
include the ordinary K\"{a}hler as a special case.  They come from 
transforming both
the complex and symplectic structures which occur in the ordinary
K\"{a}hler case by the $B$-field $B$.

\begin{equation}\label{e:JJ}
\cJ=\left(
\begin{array}{cc}
J & 0  \\
BJ + J^{\vee}B & -J^{\vee} 
\end{array}
\right)
\end{equation}

\noindent
and

\begin{equation}\label{e:Jom}
\cJ'=\left(
\begin{array}{cc}
\om^{-1}B & -\om^{-1}  \\
\om + B \om^{-1}B & -B\om^{-1}  
\end{array}
\right)
\end{equation}

\noindent
where $\om J = -J^{\vee}\om$.  The ordinary K\"{a}hler case of course 
comes about from setting $B$ to zero.
 
\end{examples}

\

\bigskip

\section{T-duality}\label{s:T-Duality}

Our main goal is to extend  the usual T-Duality transformation of
geometric structures on families of tori in a way that will allow us 
to incorporate generalized (almost) complex structures.

\subsection{T-duality in all directions}\label{ss:duality} 

In its simplest form, T-Duality exchanges geometric data on a torus
$\boldsymbol{T} \cong (S^{1})^{\times n}$ with geometric data on the
dual torus $\boldsymbol{T}^{\vee}$. For instance if the torus
$\boldsymbol{T}$ is a complex manifold, then the dual
torus is also naturally a complex manifold. This
immediately generalizes to translation invariant (hence integrable)
generalized complex structures on $\boldsymbol{T}$.

Indeed, choose a realization of $\boldsymbol{T}$ as a
quotient $\boldsymbol{T} = V/\Lambda$ of a real $n$-dimensional vector
space $V$ by a sub-lattice ${\mathbb Z}^{n} \cong \Lambda \subseteq
V$. Then specifying a translation invariant generalized complex
structure on $\boldsymbol{T}$ is equivalent to specifying a constant
generalized complex structure $\cJ \in GL(V\oplus V^{\vee})$ on the
vector space $V$. Now the dual torus $\boldsymbol{T}^{\vee}$ has a
natural realization as the quotient $\boldsymbol{T}^{\vee} =
V^{\vee}/\op{Hom}(\Lambda,{\mathbb Z})$. Thus, in order to describe
the $T$-dual generalized complex structure on $\boldsymbol{T}^{\vee}$
it suffices to produce a constant generalized complex structure on
$V^{\vee}$. This can be done in a simple way:  
Let $\tau:V\oplus V^{\vee}\to V^{\vee}\oplus V$ be the transposition of the
two summands. Using the natural identification of $V^{\vee\vee}$
with $V$, we can also view $\tau$ as an isomorphism between
$V\oplus V^{\vee}$ and $V^{\vee} \oplus V^{\vee\vee}$. We will continue to 
denote
by $\tau$ the induced isomorphism $V_{\bC}\oplus V^{\vee}_{\bC}\to
V^{\vee}_{\bC}\oplus V^{\vee\vee}_{\bC}\cong V^{\vee}_{\bC}\oplus V_{\bC}$.  
With this notation one has the following proposition.

\begin{prop} \cite{BBB}
The isomorphism $\tau$ induces a bijection between generalized
complex structures on $V$ and generalized
complex structures on $V^{\vee}$.  If $E$ corresponds to
$\cJ\in\Aut_{\bR}(V\oplus V^{\vee})$,
then $\tau(E)$ corresponds to $\tau\circ\cJ\circ\tau^{-1}$.
\end{prop}

\begin{rem}
Below, we will see that the transformation of the spinor 
representatives is 
a Fourier-Mukai type of transformation.  The precise form of this 
transformation is given in equation \ref{e:vbspinortransform}.  Notice 
that this proposition applies equally to generalized complex structures 
on the vector space $V$ and to constant generalized complex structures 
(which are automatically integrable) on $V$ thought of as a manifold.  These 
in turn give generalized complex structures on tori which are quotients 
of the vector space.
\end{rem}
\noindent
We also have the following remark from \cite{BBB}:

\begin{rem}\label{r:dualityBbeta}
Suppose that $E$ is a generalized complex structure on a real vector
space $V$ and $E'$ is the $B$-field transform of $E$ defined by
$B\in\Wedge^2 V^{\vee}$.  Then, obviously, $\tau(E')$ is the
$\be$-field transform of $\tau(E)$, defined by the same $B\in\Wedge^2
V^{\vee}$ (but viewed now as a bi-vector on $V^{\vee}$). Thus, the
operation $\tau$ interchanges $B$- and $\be$-field transforms.
\end{rem}

The relationship from this last remark was exploited in \cite{K} to
produce an interesting conjectural relationship to non-commutative
geometry.

\subsection{More general T-duality}\label{ss:genduality} 

It has been known for some time that the previous example of
T-duality generalizes immediately to a whole family of T-duality
transformations.  This can be found for example \cite{K2} and 
the references therin.  More recently Tang and Weinstein \cite{WeiTan}
applied this observation to Dirac structures to investigate the group
of Morita equivalences of real non-commutative tori.  

By analogy with the Tang-Weinstein construction we note
that if $V = \bigoplus_{i=1}^{m} V_{i}$ and $W = \bigoplus_{i=1}^{m}
W_{i}$, where each $W_{i}$ equals either $V_{i}$ or ${V_{i}}^{\vee}$,
then the obvious isomorphism $\tau$ from $V \oplus V^{\vee}$ to $W \oplus
W^{\vee}$ intertwines the canonical quadratic forms and hence it
similarly gives a bijection between generalized complex structures on
$V$ with those on $W$.  Notice that these transformations are all real
and so there is no problem with the transversality condition.  In general, 
one could also consider as
duality transformations, isometries $\tau$, from $V_{\bC} \oplus
V_{\bC}^{\vee}$ to $W_{\bC} \oplus W_{\bC}^{\vee}$ such that
$\tau\circ\cJ\circ\tau^{-1}$ is a generalized complex structure on $W$
for all (or a family of) generalized complex structures $\cJ$ on $V$.  
A special case of this duality can easily be seen to be the right starting 
point in generalizing the symplectic/complex correspondence in \cite{Pol}.  
To see this, let $M$ be a real manifold with trivial tangent bundle, $X$ a 
real torus with its normal group structure and $V$ the tangent space to $X$ 
at the identity, thought of as a trivial bundle on $M$.  Let $\widehat{X}$ 
be the dual torus to $X$.  Then 
$T_{M \times X} \cong \pi^*(T_M \oplus V)$, and 
$T_{M \times \widehat{X}} \cong \pih^*(T_M \oplus V^{\vee})$, so for 
any isomorphism $L : T_M \to V$ we have that  

\begin{equation}
\pi^*\left(\begin{array}{cc}
0 & L   \\
-L^{-1} & 0  
\end{array}\right)
\end{equation}

\noindent
is a complex structure on $M \times X$ and 

\begin{equation}
\pih^*\left(\begin{array}{cc}
0 & L   \\
-L^{\vee} & 0  
\end{array}\right)
\end{equation}
\noindent
is a symplectic structure on $M \times \widehat{X}$.  Before pulling back,
these structures, thought of as generalized complex structures as in example 
\ref{d:complsympl} on $V \oplus T_M$ and $V^{\vee} \oplus T_M$, are related by the obvious map
 
\[
V \oplus T_M \oplus V^{\vee} \oplus T_M^{\vee} \to 
V^{\vee} \oplus T_M \oplus V \oplus T_M^{\vee}.
\]

\medskip

\section{Mirror partners of generalized 
almost complex and generalized almost K\"{a}hler structures}\label{s:vectmirror}

In this section we consider a manifold $M$ equipped with a real 
vector bundle $V$ 
where the rank of $V$ equals the dimension of $M$.  For any 
connection $\nabla$ on 
$V$ we show how to build generalized almost complex structures on 
$X = tot(V)$ in terms of data on the base manifold $M$.  We show that there is 
a bijective correspondence between generalized almost complex structures 
built in this way on $X$ and generalized almost complex structures of the same 
type on $\Xwh = tot(V^{\vee})$ built using $\nabla ^{\vee}$.

Let $X$ be the total space of any vector bundle $V$ over a manifold $M$.  
Then we have the exact tangent sequence 

\begin{equation}\label{e:tanseq}
0 \rar{} \pi^*V \rar{j} T_X \rar{d\pi}
\pi^*T_M \rar{} 0
\end{equation}
A connection on the bundle $V$ is by definition a map of sheaves

\[
V \rar{\nabla} V \otimes T_M^{\vee} 
\]

\noindent
satisfying $\nabla(f\sigma) = \sigma \otimes df + f \nabla(\sigma)$ for 
all local 
sections $f$ of $C^\infty_M$ and 
$\sigma$ of $V$.  We can use any such connection to give a splitting 
of the above tangent sequence.  Namely, 
let 
\[
\pi^*\nabla: \pi^* V \rightarrow \pi^* V \otimes T_X^{\vee}
\] 
be the pullback 
of $\nabla$ and let $S$ be the tautological global section of $\pi^*V$ 
on $X$.  
Then $D = (\pi^*\nabla)(S)$ provides a map of vector bundles 
$\pi^*V \leftarrow T_X$. 
Now its easy to see 
that this map is a splitting of \ref{e:tanseq}.  Indeed, given a 
local frame $\lbrace e_i \rbrace$ of $V$ over an open set $U \subseteq M$,  
define smooth functions $\xi_i$ on $\pi^{-1}(U)$ by 
$\xi_i (a_j e_j(m)) = a_i$ for each $m$ in $M$.  Together with the functions
$x_i \circ \pi$, for $\lbrace x_i \rbrace$ coordinates on $U \subseteq M$, 
these form a 
coordinate system in $\pi^{-1}(U)$ in which we have
$j(e_i) = \pdxii$.  In these coordinates we have that on $\pi^{-1}(U)$, 

\[
S=\xi_i \pi^{-1}e_i
\] 

\noindent
and so if we define $D$ by
 
\begin{equation}\label{eqn:Ddef}
D =(\pi^*\nabla) (S) 
= \pi^{-1}e_i \otimes d\xi_i + \xi_i \pi^{-1}e_j \otimes\pi^*A_{ji}.
\end{equation}

\noindent  
then since $\pi^*A_{ji}$ annihilates the image of $j$ we have that 
\[
D(j(\pi^{-1}e_k)) = (\pi^{-1}e_i)(d\xi_i j(\pi^{-1}e_k)) = \pi^{-1}e_k
\] 
and so $D \circ j$ is the identity.  We will write this splitting on $X$ as

\[
\xymatrix@1{0 \ar[r] & \ar[r] \pi^*V \ar[r]^-{j}  & T_X \ar@/^0.5pc/[l]^-{D} 
\ar[r]^-{d\pi} &
\ar[r]^-{} \ar@/^0.5pc/[l]^-{\al} \ar[r] \pi^*T_M \ar[r] & 0},
\]

\sbr

\noindent
Consider the isomorphism 

\begin{equation}\label{e:F}
F: T_X \oplus T_X^{\vee} \rar{} \pi^*V \oplus \pi^*T_M \oplus \pi^*V^{\vee} \oplus \pi^*T_M^{\vee}, \ \ \ \ \ \    
F=\left(
\begin{array}{cc}
D & 0  \\
d\pi & 0 \\
0 & j^{\vee} \\
0 & \alpha^{\vee}
\end{array}
\right).
\end{equation}

\noindent
with inverse

\begin{equation}\label{e:Fi}
F^{-1}:\pi^*V \oplus \pi^*T_M \oplus \pi^*V^{\vee} \oplus \pi^*T_M^{\vee} \rar{} T_X \oplus T_X^{\vee}, \ \ \ \ \ \
F^{-1}=\left(
\begin{array}{cccc}
j & \alpha & 0 & 0  \\
0 & 0 & D^{\vee} & (d\pi)^{\vee}  
\end{array}
\right).
\end{equation}

\noindent
These maps intertwine the obvious quadratic forms and therefore if 
$\underline{\cJ}$ is a generalized almost complex structure on $V \oplus T_M$ 
then $\cJ = F^{-1} (\pi^*\underline{\cJ}) F$ is a generalized almost 
complex structure on $X$.

\

\begin{defin} \label{def:nablalift}
If $\nabla$ is any connection on $V$ then we define a 
{\em $\nabla$-lifted generalized almost complex structure} to be 
a generalized almost complex structure on $X = tot(V)$ which can be 
expressed as 
$\cJ = F^{-1} (\pi^*\underline{\cJ}) F$ where $\underline{\cJ}$ is a 
generalized almost complex structure on $X$ and $F$ depends on $\nabla$ 
as explained above.
\end{defin}

\


Now using the dual connection $\nabla^{\vee}$, we may split the 
sequence tangent sequence of $\Xwh$ as

\[
\xymatrix@1{0 \ar[r] & \ar[r] \pih^*V^{\vee} \ar[r]^-{\jh}  & T_{\Xwh} \ar@/^0.5pc/[l]^-{\Dwh} 
\ar[r]^-{d\pih} &
\ar[r]^-{} \ar@/^0.5pc/[l]^-{\alh} \ar[r] \pih^*T_M \ar[r] & 0},
\]

\noindent
Of course we will also need the maps

\begin{equation}\label{e:Fwh}
\Fwh:T_{\Xwh} \oplus T_{\Xwh}^{\vee} \rar{} \pi^*V^{\vee} \oplus \pih^*T_M \oplus \pih^*V \oplus \pih^*T_M^{\vee}, \ \ \ \ \ \
\Fwh=\left(
\begin{array}{cc}
\Dwh & 0  \\
d\pih & 0 \\
0 & {\jh}^{\vee} \\
0 & \alh^{\vee}
\end{array}
\right)
\end{equation}

\noindent
with inverse

\begin{equation}\label{e:Fiwh}
\Fiwh:\pi^*V^{\vee} \oplus \pi^*T_M \oplus \pi^*V \oplus \pih^*T_M^{\vee} \rar{} T_{\Xwh} \oplus T_{\Xwh}^{\vee}, \ \ \ \ \ \
\Fiwh=\left(
\begin{array}{cccc}
\jh & \alh & 0 & 0  \\
0 & 0 & \Dwh^{\vee} & (d\pih)^{\vee}  
\end{array}
\right)
\end{equation}

\noindent
Now if we take any 
\[
\cK \in GL(V \oplus T_M \oplus V^{\vee} \oplus T_M^{\vee})
\]
we can apply the duality transformation along the fibers to get
\[
\cKwh \in GL(V^{\vee} \oplus T_M \oplus V \oplus T_M^{\vee}).
\]

\noindent
Clearly this transformation intertwines the quadratic forms and so 
$\cK$ is a generalized almost complex structure on $V \oplus T_M$ if and only 
if $\cKwh$ is a generalized almost complex structure on $V^{\vee} \oplus T_M$.  Therefore $\cJ = F^{-1} (\pi^*\cK) F$  
is a generalized almost complex structure on $X$ 
if and only  $\cJwh = \Fiwh ({\pih}^*\cKwh) \Fwh$ is a 
generalized almost complex structure on $\Xwh$.  At this point we 
will impose an extra constraint on these structures.  

\

\begin{defin} \label{def:nablaadapted}
A $\nabla$-lifted generalized almost complex structure 
$\cJ = F^{-1} (\pi^*\underline{\cJ}) F$ will be called 
{\em adapted} if 
\[
\cK (V \oplus V^{\vee}) = T_M \oplus T_M^{\vee}
\]
\end{defin}

\

\noindent
We will assume that $\cJ$ is an adapted, $\nabla$-lifted generalized 
almost complex structure from now on.

\begin{rem}
It is clear from the construction above that $\cJ$ is adapted if and 
only if  $\cJwh$ is.
\end{rem}

Finally, let us record the explicit formulas for the operators 
$\cK$, $\cKwh$, $\cJ$ and $\cJwh$.  The adapted condition together with the 
fact that ${\cK}^2 = -1$ and that $\cK$ preserves the quadratic form ensure 
that it is of the form  
\begin{equation}\label{e:renametoK}
\cK=\left(
\begin{array}{cccc}
0             & \cJ_{12}       & 0                  & \cJ_{22} \\
\cJ_{13}           & 0         & -\cJ_{22}^{\vee}        & 0 \\
0             & \cJ_{31}       & 0                  & -\cJ_{13}^{\vee} \\
-\cJ_{31}^{\vee}   & 0         & -\cJ_{12}^{\vee}        & 0
\end{array}
\right), \ \ \ \ \cK \in GL(V \oplus T_M \oplus V^{\vee} \oplus T_M^{\vee})
\end{equation}

\noindent
subject to 

\begin{eqnarray}
\label{e:K1} &&  \cJ_{12} \cJ_{13}          - \cJ_{22} \cJ_{31}^{\vee} = -1; \\
\label{e:K2} &&  \cJ_{12} \cJ_{22}^{\vee}   +\cJ_{22} \cJ_{12}^{\vee} = 0; \\
\label{e:K3} &&  \cJ_{13} \cJ_{12}          - \cJ_{22}^{\vee} \cJ_{31} = -1; \\
\label{e:K4} &&  \cJ_{13} \cJ_{22}          + \cJ_{22}^{\vee} \cJ_{13}^{\vee} = 0; \\
\label{e:K5} &&  \cJ_{31} \cJ_{13}          + \cJ_{13}^{\vee} \cJ_{31}^{\vee}=0; \\
\label{e:K6} &&  \cJ_{31}^{\vee} \cJ_{12}   + \cJ_{12}^{\vee} \cJ_{31}=0. 
\end{eqnarray}

\noindent
With this notation we have

\begin{equation}\label{e:K_check}
\cKwh=\left(
\begin{array}{cccc}

0              & \cJ_{31}     & 0               & -\cJ_{13}^{\vee}   \\
-\cJ_{22}^{\vee}    & 0       & \cJ_{13}             & 0             \\
0              & \cJ_{12}     & 0               & \cJ_{22}           \\
-\cJ_{12}^{\vee}    & 0       & -\cJ_{31}^{\vee}     & 0

\end{array}
\right), \ \ \ \ \cKwh \in GL(V^{\vee} \oplus T_M \oplus V \oplus T_M^{\vee})
\end{equation}

\noindent
and so

\begin{equation}\label{e:J}
\cJ=\left(
\begin{array}{cc}
j(\pi^*\cJ_{12}) (d\pi) + \alpha (\pi^*\cJ_{13}) D & 
j (\pi^*\cJ_{22}) \alpha^{\vee} - \alpha (\pi^*\cJ_{22}^{\vee}) j^{\vee}   \\
D^{\vee} (\pi^*\cJ_{31}) (d\pi) -(d\pi)^{\vee} (\pi^*\cJ_{31}^{\vee}) D & 
-D^{\vee} (\pi^*\cJ_{13}^{\vee}) \alpha^{\vee} -(d\pi)^{\vee} (\pi^*\cJ_{12}^{\vee}) j^{\vee}   
\end{array}
\right).
\end{equation}

\noindent
and

\begin{equation}\label{e:cJ}
\cJwh=\left(
\begin{array}{cc}

\jh({\pih}^*\cJ_{31})(d\pih) -  
\alh(\pih^*\cJ_{22}^{\vee})\Dwh & 

-\jh({\pih}^*\cJ_{13}^{\vee})\alh^{\vee} 
+ \alh({\pih}^*\cJ_{13})\jh^{\vee}   \\

\Dwh^{\vee}({\pih}^*\cJ_{12})(d\pih) 
-(d\pih)^{\vee}({\pi}^*\cJ_{12}^{\vee}) \Dwh & 

\Dwh^{\vee}({\pih}^*\cJ_{22})\alh^{\vee} 
-(d\pih)^{\vee}({\pih}^*\cJ_{31}^{\vee})\jh^{\vee}   

\end{array}
\right).
\end{equation}

\begin{rem}\label{rem:Jswitch}
Notice that the mirror symmetry transformation ``exchanges'' 
$\cJ_{12}$ with $\cJ_{31}$ and $\cJ_{22}$ with $-\cJ_{13}^{\vee}$.  
\end{rem}
We have written down the bijective correspondence between $\nabla$-lifted, 
adapted, generalized almost complex structures on $X$ and 
$\nabla^{\vee}$-lifted, adapted, generalized almost complex structures on 
$\Xwh$.  We will show below, in the case that $\nabla$ is flat, 
that $\cJ$ is integrable if and only if $\cJwh$ is.

\subsection{Associated almost Dirac structures}\label{ss:Dirac}

For each of the generalized complex structures on $X$ that we consider, 
there is a natural almost Dirac structure that appears on the base manifold 
$M$.  It does not depend on the connection used to split the tangent sequence 
of $X \to M$.  An almost Dirac structure on $M$ is just \cite{Cou} a maximally 
isotropic sub-bundle of $T_M \oplus T_M^{\vee}$.   
Now the isomorphism $\cK$,  given in equation (\ref{e:renametoK}), 
preserves the quadratic form and when restricted to $V \oplus V^{\vee}$, gives 
an isomorphism $V \oplus V^{\vee} \to T_M \oplus T_M^{\vee}$ which 
preserves the obvious quadratic forms.  Hence the image of $V$ is a 
maximally isotropic subspace of $T_M \oplus T_M^{\vee}$.  In other words 
it is an almost Dirac structure on $M$ which we will call $\Delta$, where

\[   
\Delta = \cK(V) = \cJ_{13}(V) -\cJ_{31}^{\vee}(V) =\cKwh(V).
\]

\begin{examples}
\
\begin{enumerate}
\item
Suppose we use our method to construct an 
almost complex structure on $X = tot(V)$ out of some arbitrary 
connection on $V$.  Then we necessarily have that $\cJ_{13}$ is an isomorphism and 
$\cJ_{31} = 0$.  Hence $\Delta = T_M$.
\item 
If instead we put an almost symplectic structure on $X = tot(V)$ then 
$\Delta = T_M^{\vee}$.
\end{enumerate} 
\end{examples}

\

Notice that the almost Dirac structure 

\[
\widehat{\Delta} = \cKwh (V^{\vee}) 
= -\cJ_{22}(V^{\vee}) -\cJ_{12}^{\vee}(V^{\vee}) =  \cK(V^{\vee})
\]
arising from the mirror generalized almost complex structure is always 
transverse to $\Delta$.  Hence we always get a pair

\[   
\Delta \oplus \widehat{\Delta}= T_M \oplus T_M^{\vee}
\]

\noindent
of complementary almost Dirac structures.  Later we will return 
to these structures and study their 
integrability and the existence of flat connections on them.

\subsection{Mirror symmetry for generalized almost K\"{a}hler
  manifolds}\label{ss:genkah} 

\

In this section, we study the case of a pair of $\nabla$-lifted, adapted 
generalized almost complex structures on the total space of a 
vector bundle which 
form a generalized almost K\"{a}hler structure as described in 
\ref{ss:kahlerdef}.  Under these
conditions, we write down the mirror transformation rule that allows us 
to relate the
generalized almost K\"{a}hler metric $G$ on $X$ and the mirror generalized
almost K\"{a}hler metric $\Gwh$ on $\Xwh$.  We observe that in
general, the local transformation rules for the pair $(g,b)$ that exist in 
the physics literature, continue to hold in this setting, even though here
neither $\cJ$ nor $\cJ'$ needs to be a B-field transform of a generalized 
complex structure of complex type.  

First of all notice that the mirror transform of a generalized 
almost K\"{a}hler pair
$(\cJ, \cJ')$ is also generalized almost K\"{a}hler.  Indeed, if we let
$\cJ = F^{-1}(\pi^*\cK)F$ and $\cJ' = F^{-1}(\pi^*\cK')F$, then 
$\cJ$ and $\cJ'$ commute
if and only if $\cK$ and $\cK'$ commute.  This, in turn, is equivalent to 
$\cKwh$ and $\cKpwh$ commuting which happens if and only if 
$\cJwh = \Fiwh(\pih^*\cKwh)\Fwh$ and 
$\cJpwh = \Fiwh(\pih^*\cKpwh)\Fwh$ commute.
Similarly, $G = -F^{-1} \pi^*(\cK \cK') F$ is positive definite if and only 
if $-\cK \cK'$ is.  This is equivalent to $-\cKwh \cKpwh$ being 
positive definite, which happens if and only if 
$\Gwh = -\Fiwh \pih^*(\cKwh \cKpwh) \Fwh$ is positive definite.  
By our assumptions on  $\cJ$ and $\cJ'$ we
may write $\underline{G} = -\cK \cK'$ as

\begin{equation}\label{e:Kahler1}
G:V \oplus T_M \oplus V^{\vee} \oplus T_M^{\vee} \to V \oplus T_M \oplus V^{\vee} \oplus T_M^{\vee}, \ \ \ \ \ \
\underline{G}= -\cK\cK'=\left(
\begin{array}{cccc}
G_{11}                & 0         & G_{21}             & 0         \\
0                  & G_{14}       & 0               & G_{24}           \\
G_{31}                & 0         & G_{11}^{\vee}      &  0\\
0                  & G_{34}       & 0               & G_{14}^{\vee}
\end{array}
\right)
\end{equation}
\noindent
where, $G_{21} = G_{21}^{\vee}$,$G_{24} = G_{24}^{\vee}$,  $G_{31} = G_{31}^{\vee}$, 
$G_{34} = G_{34}^{\vee}$.  Finally, using the fact that this matrix 
squares to the identity, we get:

\[
G_{11} = -\cJ_{12} \cJ'_{12} + \cJ_{22} {\cJ'_{31}}^{\vee}
\]
\[
G_{21} = \cJ_{12} {\cJ'_{21}}^{\vee} + \cJ_{22} {\cJ'_{11}}^{\vee}
\]
\[
G_{14} = -\cJ_{13} {\cJ'_{11}}^{\vee} + \cJ_{22}^{\vee} \cJ'_{31}
\]
\[
G_{24} = -\cJ_{13} \cJ'_{21} -\cJ_{22}^{\vee} {\cJ'_{12}}^{\vee}
\]
\[
G_{31} = -\cJ_{31} \cJ'_{12} - \cJ_{13}^{\vee} {\cJ'_{31}}^{\vee}
\]
\[
G_{34} = \cJ_{31}^{\vee} \cJ'_{11} +\cJ_{12}^{\vee} \cJ'_{31}.
\]
\noindent
Therefore $\underline{G}' = -\cKwh \cKpwh$ comes out to be

\begin{equation}\label{e:Kahler2}
\underline{G}' = -\cKwh  \cKpwh =\left(
\begin{array}{cccc}
G_{11}^{\vee}         & 0         & G_{31}             & 0         \\
0                  & G_{14}       & 0               & G_{24}           \\
G_{21}                & 0         & G_{11}             &  0\\
0                  & G_{34}       & 0               & G_{14}^{\vee}
\end{array}
\right).
\end{equation}

\begin{rem}
Notice that the mirror symmetry transformation ``exchanges'' 
$G_{11}$ with $G_{11}^{\vee}$, $G_{21}$ with $G_{31}$, 
and ``preserves'' $G_{14}$ and $G_{34}$.
\end{rem}
\noindent
Now writing $G$ in terms of $g$ and $b$, (\cite{Gua})

\begin{equation}\label{e:the_form2}
G=\left(
\begin{array}{cc}

-g^{-1} b          & g^{-1}         \\
g - bg^{-1}b       & bg^{-1}               

\end{array}
\right)
\end{equation}

\noindent
and similarly writing $\Gwh$ in terms of $\gh$ and $\bh$ we can easily 
manipulate the resulting equations to yield the following formulas for 
the metrics and B-fields in terms of the vector bundle maps $G_{ij}$ on 
the base manifold.
\[
g = D^{\vee} \pi^*G_{21}^{-1} D + (d\pi)^{\vee} \pi^*G_{24}^{-1} d\pi
\]
\[
b = D^{\vee} \pi^*(G_{11}^{\vee} G_{21}^{-1}) D
     + {(d\pi)}^{\vee} \pi^*(G_{14}^{\vee} G_{24}^{-1}) d\pi
\]
\[
\gh = {\Dwh}^{\vee} \pih^*G_{31}^{-1} {\Dwh} 
+ {(d\pih)}^{\vee} \pih^* G_{24}^{-1} d \pih
\]
\[
\bh =  {\Dwh}^{\vee} \pih^*(G_{11} G_{31}^{-1}) {\Dwh} 
       +  {(d\pih)}^{\vee} \pih^*(G_{14}^{\vee} G_{24}^{-1}) d \pih
\]

Notice that our assumptions on the compatibility of the generalized complex 
structures, and the foliation and transverse vector bundle, imply that 
the metric $g$ and B-field $b$ do not mix the horizontal 
and vertical directions.

Now if we chose local vertical coordinates adapted to the 
flat connection, $y^\al$ on $X$ 
and ${\yh}^\al$ on $\Xwh$ and $x^i$ on the base then the above 
just means that locally we have

\[
g = g_{ij}(x) dx^i dx^j + h_{\al\beta}(x) dy^\al dy^\beta
\]
\[
b = b_{ij}(x) dx^i dx^j + B_{\al\beta}(x) dy^\al dy^\beta
\]
\[
\gh = g_{ij}(x) dx^i dx^j + {\hh}_{\al\beta}(x) d{\yh}^\al d{\yh}^\beta
\]
\[
\bh = b_{ij}(x) dx^i dx^j + {\Bwh}_{\al\beta}(x) d{\yh}^\al d{\yh}^\beta
\]
where of course, $x^i$ means $x^i \circ \pi$ on $X$ and $x^i \circ \pih$ 
on $\Xwh$.

Then the Buscher transformation rules \cite{Bus1, Bus2}
(we used \cite{Allessandro} as a reference)
\[
(h + B) \hh (h-B) = h
\ \ \ \ \ \
and
\ \ \ \ \ \
(h+B) \Bwh (h-B) = -B
\]
are verified from the easily checked identities
\[
(G_{21}^{-1} + G_{11}^{\vee} G_{21}^{-1})G_{31}^{-1}
(G_{21}^{-1} - G_{11}^{\vee}G_{21}^{-1}) = G_{21}^{-1}
\]
and
\[
(G_{21}^{-1} + G_{11}^{\vee} G_{21}^{-1})G_{11} {G_{31}}^{-1}
(G_{21}^{-1} - G_{11}^{\vee} G_{21}^{-1}) = -G_{11}^{\vee} {G_{21}}^{-1}
\]
respectively.

We now work out the transformation rules relating the two almost
complex structures, $J_+$, $J_-$, and their mirror partners $\Jwh_+$
and $\Jwh_-$.  We have 
\[
J_+ = \cJ_1 + \cJ_2 (g + b)
\ \ \ \ \ \
{\mbox and}
\ \ \ \ \ \
J_- = \cJ_1 + \cJ_2 (b - g).
\]  
By combining the results above we can easily compute that
\[
J_+ = j(\pi^*(\cJ_{12} + \cJ_{22}(G_{14}^{\vee} + 1)G_{24}^{-1}))d\pi
+\al(\pi^*(\cJ_{13} - \cJ_{22}^{\vee} (G_{11}^{\vee} + 1) G_{21}^{-1}))D
\]
{\mbox and}
\[
J_- = j(\pi^*(\cJ_{12} + \cJ_{22}(G_{14}^{\vee} - 1)G_{24}^{-1}))d\pi
+\al(\pi^*(\cJ_{13} - \cJ_{22}^{\vee} (G_{11}^{\vee} - 1) G_{21}^{-1}))D.
\]
Hence

\[
\Jwh_+ = \jh(\pih^*(\cJ_{31} -\cJ_{13}^{\vee}(G_{14}^{\vee} + 1)G_{24}^{-1}))d\pih
+\alh(\pih^*(-\cJ_{22}^{\vee} + \cJ_{13} (G_{11} + 1) G_{31}^{-1}))\Dwh
\]
and

\[
\Jwh_- = \jh(\pih^*(\cJ_{31} - \cJ_{13}^{\vee}(G_{14}^{\vee} - 1)G_{24}^{-1}))d\pih
+\alh(\pih^*(-\cJ_{22}^{\vee} + \cJ_{13} (G_{11} - 1) G_{31}^{-1}))\Dwh.
\]

\section{Branes}\label{s:Branes}

In this section, we give some ideas of how one can transfer branes 
\cite{Gua, K} from a generalized almost complex manifold to 
its mirror partner.  
We will present in detail only a very restricted case.  This construction 
closely parallels that in 
\cite{Pioli, Leung}.  Consider the following definition from \cite{K} which 
also appears in a more general form in \cite{Gua}.

\begin{defin} \label{def:branedef} \cite{K}
Let $(X, \cJ)$ be a generalized (almost) complex manifold.  Consider triples 
\[
(Y,\cL,\nabla_{\cL})
\]
 where $f: Y \hookrightarrow X$ is a 
sub-manifold of $X$, $\cL$ 
is a Hermitian line bundle on $Y$, and
$\nabla_{\cL}$ is a connection on $\cL$.  Such a triple is said 
to be a {\em generalized complex brane} if the bundle  
\[
\{(v,\alpha) \in T_Y \oplus (T_{X}^{\vee} |_{Y})\ \   |  \ \  
f^* \circ (df)^{\vee} \alpha = \myiota_v\cF \}
\]
is preserved by the restriction of $\cJ$ to $Y$, where $\cF$ is the curvature 
two-form of $\nabla_{\cL}$. 
\end{defin}

We studied some special cases \cite{Gua} of these branes 
in \cite{BBB} under the name of 
generalized Lagrangian sub-manifolds and 
found some interesting relationships to sub-manifolds of $X$ which inherit 
generalized complex structures (which we call generalized 
complex sub-manifolds).

Suppose that $M$ is an $n$-manifold, $V$ is a rank $n$ vector bundle on $M$, 
$\nabla$ is a connection on $V$, $X$ is the total space of $V$, 
$\Xwh$ is the total space of $V^{\vee}$  
and $\cJ$ is an adapted, $\nabla$-lifted (see section \ref{s:vectmirror}) 
generalized almost complex structure on $V$.  Let $S$ be a sub-manifold of $M$, 
$W \subseteq V |_{S}$ a sub-bundle, $Y$ the total space of $W$, and $\Ywh$ 
the total space of the sub-bundle $Ann(W) \subseteq V^{\vee} |_{S}$.  Then 
we propose that the relationship between $Y$ and $\Ywh$ is a special case of 
a potential generalization of the relationship 
between $A-$ cycles and $B-$ cycles in mirror symmetry 
(see e.g. \cite{Leung} and references therein).  We justify 
this with the following lemma.

\begin{lem}\label{lem:branetransform}
Under the conditions of the preceding paragraph, 
the triple $(Y, \mathbb{C} \otimes C^{\infty}_Y, d)$ is a generalized 
complex brane of $(X, \cJ)$ 
if and only if the triple 
$(\Ywh, \mathbb{C} \otimes C^{\infty}_{\Ywh}, d)$ is a 
generalized complex brane of $(\Xwh, \cJwh)$.
\end{lem}

\noindent
{\bf Proof.}

\

In this proof, we will be using the notation of section \ref{s:vectmirror}.
We need to show that 
\begin{equation}\label{e:thing1}
\cJ(T_Y \oplus \Ann(T_Y)) = T_Y \oplus \Ann(T_Y)
\end{equation}
if and only if 
\begin{equation}\label{e:thing2}
\cJwh(T_{\Ywh} \oplus \Ann(T_{\Ywh})) =
 T_{\Ywh} \oplus \Ann(T_{\Ywh}), 
\end{equation}
where it is to be understood that we are restricting $\cJ$ to $Y$ and 
$\cJwh$ to $\Ywh$.  Observe that when understood as bundles on $Y$, we have
\[
T_Y = j(\pi^* W) \oplus \al(\pi^*T_S), \ \ 
\ \ \Ann(T_Y) = D^{\vee}(\pi^*(\Ann(W))) \oplus (d\pi)^{\vee}(\pi^*(\Ann(T_S)))
\] 
and when understood as bundles on $\Ywh$, we have
\[
T_{\Ywh} = \jh(\pih^* (\Ann(W))) \oplus \alh(\pih^*T_S), \ \ 
\ \ \Ann(T_Y) = \Dwh^{\vee}(\pih^*(W)) \oplus (d\pih)^{\vee}(\pih^*(\Ann(T_S))).\]
From this perspective it is clear that both \ref{e:thing1} and \ref{e:thing2} 
are both equivalent simply to the conditions (understanding that 
$\cK$ is restricted to $S$)

\begin{eqnarray}
\label{e':1} && \cJ_{13}(W) \subseteq T_S \\
\label{e':2} && \cJ_{12}(T_S) \subseteq W\\
\label{e':3} && \cJ_{22}(\Ann(T_S)) \subseteq W\\
\label{e':4} && \cJ_{31}(T_S) \subseteq \Ann(W).
\end{eqnarray}

\noindent
and therefore we are done.  A more general treatment will involve replacing
$W$ by an affine sub-bundle which will result in a non-trivial line bundle 
on the mirror side.  A more general story will be the 
result of upcoming work.  The extension of the results in this section 
to the case of torus bundles will be left to the reader, and should be clear 
upon reading section \ref{s:tbundle}.   We hope that 
in a suitable extended version of the homological mirror symmetry 
conjecture \cite{Ko}, generalized complex manifolds would be 
assigned categories in a natural way, and branes would be related to 
objects in these categories.    

\begin{rem}
On a torus bundle $Z \to M$, where $Z$ is an orientable compact manifold, 
it is plausible that the correspondence which we are describing here, when 
thought of as a correspondence between homology classes on $Z$ to homology 
classes on the dual torus bundle $\Zwh \to M$ agrees, upon using Poincar\'{e} 
Duality, with the correspondence in cohomology given in 
section~\ref{s:tbundle}.  
\end{rem}

\section{The mirror transformation on spinors and the Fourier
  transform}\label{s:fourier} 

In this section we study a map from certain complex valued
differential forms on the total space of a vector bundle to complex
valued differential forms on the total space of the dual vector
bundle.  We show that the line sub-bundle of the bundle of
differential forms associated to an adapted, $\nabla$-lifted
generalized almost complex structure has a sub-sheaf which goes under
this correspondence to the sub-sheaf associated to the mirror
generalized almost complex structure.  The idea of using a Fourier 
transform in the context of T-duality for generalized complex 
structures has appeared in a slightly different context in both \cite{Gua} 
(based on ideas appearing in \cite{Mathai}) and also \cite{Allessandro}, 
\cite{Van} and the references therin.

Consider a vector bundle $V$ of rank $n$ on a manifold $M$.  There is 
an isomorphism 

\[
\Wedge V^{\vee} \otimes{\mathbb{C}} \rar{} 
(\Wedge V \otimes \Wedge^n V^{\vee})\otimes{\mathbb{C}}
\]
\noindent given by 
\begin{equation}\label{e:vbspinortransform}
\phi \mapsto \int (\phi \wedge \exp(\kappa))
\end{equation}
where $\kappa$ is the canonical global section of 
$(V \otimes V^{\vee})\otimes{\mathbb{C}} \subseteq 
\Wedge^2 (V \oplus V^{\vee})\otimes{\mathbb{C}}$
and 
\[\int : \Wedge (V \oplus V^{\vee})\otimes{\mathbb{C}} 
\to \Wedge V \otimes \Wedge^n V^{\vee}\otimes{\mathbb{C}}
\]
is the projection map.  Furthermore, this map decomposes into isomorphisms

\[
\Wedge^p V^{\vee}\otimes{\mathbb{C}} \to \Wedge^{n-p} V \otimes \Wedge^n V^{\vee} \otimes{\mathbb{C}}.
\] 

\noindent
and also induces a Fourier transform isomorphism, which we 
will call {\bf F.T.} 

\[
\Wedge V^{\vee} \otimes \Wedge T_M^{\vee} \otimes{\mathbb{C}} 
\rar{\mathbf{F.T.}} \Wedge V \otimes \Wedge 
T_M^{\vee} \otimes \Wedge^n V^{\vee} \otimes{\mathbb{C}}
\]

\noindent
which in turn decompose into isomorphisms

\[
(\Wedge^q T_M^{\vee} \otimes 
\Wedge^p V^{\vee} \otimes{\mathbb{C}}) 
\to (\Wedge^q T_M^{\vee} \otimes 
\Wedge^{n-p} V \otimes \Wedge^n V^{\vee}\otimes{\mathbb{C}}).
\]

\begin{lem}\label{lem:vbspin}
Let $V$ be a rank $n$ orientable vector bundle on an $n-$manifold, 
$\cK$ a generalized almost complex structure on the vector bundle 
$V \oplus T_M$ satisfying 
$\cK(V \oplus V^{\vee}) = T_M \oplus T_M^{\vee}$, and $\cKwh$ the mirror 
structure.  Then the composition of  
the map {\bf F.T.} with any trivialization of $\Wedge^n V^{\vee}$ takes 
the line bundle $\underline{L} \subseteq \Wedge(V \oplus T_M)^{\vee}\otimes{\mathbb{C}}$ 
which represents $\cK$ to the line bundle $\underline{\Lwh} 
\subseteq \Wedge(V \oplus T_M)^{\vee}\otimes {\mathbb{C}}$ which 
represents $\cKwh$. 
\end{lem}
\noindent
{\bf Proof.}

First of all notice that changing the trivialization of $\Wedge^n V$ 
multiplies the image of {\bf F.T.} 
by a non-zero function on the base manifold. This is  
an automorphism of image of the composed map along with its inclusion into 
$\Wedge(V \oplus T_M)^{\vee}\otimes {\mathbb{C}}$.  The 
$(+i)$ eigenbundle of $\cK$ is the graph of the map
\[
-\smo \cK \big\vert_{(V \oplus V^{\vee})\otimes{\mathbb{C}}} 
: (V \oplus V^{\vee})\otimes{\mathbb{C}} 
\to (T_M \oplus {T_M}^{\vee})\otimes{\mathbb{C}}.
\]
Similarly, the (+i) eigenbundle of $\cKwh$ is the graph of
\[
-\smo \cK \big\vert_{(V^{\vee} \oplus V)\otimes{\mathbb{C}}} 
: (V \oplus V^{\vee})\otimes{\mathbb{C}} 
\to (T_M \oplus {T_M}^{\vee})\otimes{\mathbb{C}}.
\]
\noindent
The sections $\phi$ of $L$ therefore satisfy

\begin{equation}\label{e:spineq1}
\myiota_{v - \smo \cJ_{13} v} \phi 
+ \smo ({\cJ_{31}}^{\vee} v) \wedge \phi =0
\end{equation}
\begin{equation}\label{e:spineq2}
\myiota_{\smo {\cJ_{22}}^{\vee} \al} \phi 
+  (\al +\smo \cJ_{12}^{\vee} \al) \wedge \phi =0 
\end{equation}
\noindent
for all sections $v$ of $V$ and $\al$ of $V^{\vee}$.

We need to show that the section $\mathbf{ F.T.}(\phi)$ of 
$\Wedge(V^{\vee} \oplus T_M)^{\vee}\otimes 
\Wedge^n V^{\vee} \otimes{\mathbb{C}}$ satisfies

\begin{eqnarray}
\label{spin:3} && \myiota_{-\smo \cJ_{13} v} \mathbf{F.T.} (\phi) 
+  (v +\smo \cJ_{31}^{\vee} v) \wedge \mathbf{F.T}(\phi) =0     \\
\label{spin:4} &&  
\myiota_{\al +\smo {\cJ_{22}}^{\vee} \al} \mathbf{F.T.}(\phi) 
+  (\smo \cJ_{12}^{\vee} \al) \wedge \mathbf{F.T.}(\phi) =0 
\end{eqnarray}
\noindent for all sections $v$ of $V$ and $\al$ of $V^{\vee}$.  These 
equations hold for the map {\bf F.T.} if and only if they hold for 
the composition of {\bf F.T.} with any trivialization.  This 
can be seen by writing {\bf F.T.} as the 
composed map followed by the 
action of ``wedging'' with a global section of $\Wedge^n V^{\vee}$.
The equations \ref{spin:3} and \ref{spin:4} will follow immediately 
from taking the Fourier transform of both sides 
of \ref{e:spineq1} and \ref{e:spineq2} and using the following lemma. 
\hfill $\Box$

\

\begin{lem} \label{lem:fouriermukailem}
For any sections $\zeta$ of $\Wedge(V \oplus T_M)^{\vee}\otimes{\mathbb{C}}$,
$v$ of $V\otimes{\mathbb{C}}$, $w$ of $T_M\otimes{\mathbb{C}}$,  
$\al$ of $V^{\vee}\otimes{\mathbb{C}}$ and $\beta$ of $T_M^{\vee}\otimes{\mathbb{C}}$ we have
\medskip

{\bf (i)}  $\mathbf{F.T.}(\myiota_{v} \zeta) = v \wedge \mathbf{F.T.}(\zeta)$
\medskip

{\bf (ii)} $\mathbf{F.T.}(\myiota_{w} \zeta) = \myiota_w \mathbf{F.T.}(\zeta)$
\medskip

{\bf (iii)} $\mathbf {F.T.}(\al \wedge \zeta) = \myiota_{\al} \mathbf{F.T.}(\zeta)$
\medskip

{\bf (iv)} $\mathbf{F.T.}(\beta \wedge \zeta) = \beta \wedge \mathbf{F.T.}(\zeta)$
\end{lem}

\noindent
{\bf Proof.}

It clearly suffices to prove this in the case that $\zeta$ is a section of 
$\Wedge^p(V \oplus T_M)^{\vee}$.  Notice 
also that $\int\ \ \circ\ \ \myiota_v = 0$ and $\myiota_v \kappa = -v$ 
for any section $v$ of $V$ and $\myiota_{\al} \kappa = \al$ for any 
section $\al$ of $V^{\vee}$.  Then we have

\begin{eqnarray*}
\mathbf{F.T.} (\myiota_{v} \zeta) \ \ 
& = & \int (\myiota_{v} \zeta) \wedge \exp(\kappa) 
= \int \myiota_{v} (\zeta \wedge \exp(\kappa)) 
- (-1)^p \int \zeta \wedge \myiota_{v} \exp(\kappa) \\
& = &- (-1)^p \int \zeta \wedge \myiota_{v} \exp(\kappa)
= - (-1)^p \int \zeta \wedge (\myiota_{v} (\kappa)) 
\wedge \exp(\kappa) \\
& = &(-1)^p \int \zeta \wedge v 
\wedge \exp(\kappa) =\int  v \wedge \zeta
\wedge \exp(\kappa) = v \wedge \int \zeta
\wedge \exp(\kappa) \\
& = &v \wedge \mathbf{F.T.} (\zeta) \\
\mathbf{F.T.} (\myiota_{w} \zeta) \ \
& = &\int (\myiota_{w} \zeta) \wedge \exp\kappa)=\int \myiota_{w} (\zeta \wedge \exp(\kappa)) 
=\myiota_{w} \int (\zeta \wedge \exp(\kappa)) \\
& = &\myiota_{w} \mathbf{F.T.}(\zeta) \\
\mathbf{F.T.} (\al \wedge \zeta) 
& = &\int  
(\al \wedge \zeta \wedge \exp(\kappa))
=(-1)^p \int  
(\zeta \wedge \al \wedge \exp(\kappa)) \\
& = &(-1)^p \int
(\zeta \wedge \myiota_{\al} \exp(\kappa))
=\myiota_{\al} \int  
(\zeta \wedge \exp(\kappa)) \\
& = &\myiota_{\al} \mathbf{F.T.}(\zeta) \\
\mathbf{F.T.} (\beta \wedge \zeta) 
& = &\int \beta \wedge \zeta \wedge \exp(\kappa)
=\beta \wedge \int \zeta \wedge \exp(\kappa) \\
& = &\beta \wedge \mathbf{F.T.}(\zeta)
\end{eqnarray*}   \hfill $\Box$

Let $M$ an $n$-manifold and $X \rar{\pi} M$ be the total space of an 
orientable vector bundle $V$ on $M$, with connection $\nabla$ and 
$\cJ$ a $\nabla$-lifted, adapted  
generalized almost complex structure on $X$.  Let 
$\Xwh \rar{\pih} M$ be the total 
space of $V^{\vee}$.   
Using $\nabla$ we may realize 
$\pi^* (\Wedge V^{\vee}
\otimes \Wedge T_M^{\vee} \otimes{\mathbb{C}})$ 
as a sub-bundle
of $\Wedge T_X^{\vee} \otimes{\mathbb{C}}$.  
Now $\cJ$ determines a spinorial line bundle 
$L \subseteq \Wedge T_X^{\vee} \otimes \mathbb{C}$ which is simply the 
image under this isomorphism of the pullback $\pi^* \underline{L}$.  Similarly,
$\cJwh$ determines a spinorial line bundle 
$\Lwh \subseteq \Wedge T_{\Xwh}^{\vee} \otimes \mathbb{C}$ isomorphic to 
$\pih^*\underline{\Lwh}$.  Therefore, interpreting the Fourier 
transform maps as isomorphisms  
$\pi_* \pi^{-1} \underline{L} \to \pih_* \pih^{-1} \underline{\Lwh}$ 
we can map certain sections of $L$ over open sets of the form 
$\pi^{-1}(U)$ to sections of $\Lwh$ over open sets of the form 
$\pih^{-1}(U)$.  We have shown the following lemma.

\begin{lem} \label{lem:fouriermukai}
Let $V$ is an orientable rank $n$ vector bundle on an $n$-manifold $M$ and 
$\cJ$ an adapted, $\nabla$-lifted generalized almost complex structure on 
$X = tot(V)$ with associated line bundle $L$.  Let the mirror generalized 
almost complex structure have associated line bundle $\Lwh$.  Then 
their are sub-sheaves, $\pi^{-1}\underline{L} \subseteq L$ and 
$\pih^{-1}\underline{\Lwh} \subseteq \Lwh$ such that if we compose the isomorphism 
\[
\pi_{*} \pi^{-1} (\Wedge T_M^{\vee} 
\otimes \Wedge V^{\vee}\otimes{\mathbb{C}} ) 
\rar{\mathbf{F.T.}} \pih_{*} \pih^{-1}(\Wedge T_M^{\vee} \otimes 
\Wedge V \otimes \Wedge^n V^{\vee} \otimes{\mathbb{C}}) 
\]

\noindent
with any trivialization of $\Wedge^n V^{\vee}$, the resulting isomorphism 
\[
\pi_*(\Wedge T_{X}^{\vee} \otimes \mathbb{C}) \supseteq 
\pi_{*} \pi^{-1} (\Wedge T_M^{\vee} 
\otimes \Wedge V^{\vee}\otimes{\mathbb{C}} ) 
\rar{} \pih_{*} \pih^{-1}(\Wedge T_M^{\vee} \otimes 
\Wedge V  \otimes{\mathbb{C}}) 
\subseteq \pih_*(\Wedge T_{\Xwh}^{\vee} \otimes \mathbb{C})
\]

\noindent
restricts to an isomorphism 
\[
\pi_*L \supseteq \pi_*\pi^{-1} \underline{L} \to \pih_{*} \pih^{-1} \underline{\Lwh}  \subseteq  \pih_*{\Lwh}
\]
\end{lem}       \hfill $\Box$

This is useful because, from the $\nabla$-lifted property, its easy to 
see that $L=\pi^{-1} \underline{L} \otimes C^{\infty}_X$ and 
$\Lwh=\pih^{-1} \underline{\Lwh} \otimes C^{\infty}_{\Xwh}$.  Therefore for $U$ small 
enough, representative spinors for $\cJ$ over $\pi^{-1}(U)$ and $\cJwh$ over 
$\pih^{-1}(U)$ exist and can be 
chosen as pullbacks of sections of $\underline{L}$ and $\underline{\Lwh}$ over $U$.  They are 
exchanged under the Fourier transform even though we have not 
written down a map between the pushforwards of $L$ and $\Lwh$.  The situation 
will be much more simple in the case of torus bundles.

\begin{rem}
It is important to remember that the geometry of $\cJ$ is not just captured by 
the abstract line bundle $L$ up to isomorphism, but rather, by $L$ together 
with its embedding into the differential forms.
\end{rem}

Understanding mirror symmetry in terms of a relationship between pure spinors 
was approached with similar techniques in \cite{Allessandro}.

\section{The question of integrability}\label{s:integrability}

The purpose of this section is to express the integrability of $\nabla$-lifted, adapted generalized almost complex structures $\cJ$ on the total space of 
vector bundles in terms of data on the base manifold $M$.  We do this only 
in the case where $\nabla$ is flat (in which case we call the structures 
$\cJ$ {\em semi-flat}.  Once we do this it will be clear that $\cJ$ 
on $X = tot(V)$ is integrable if and only if the mirror structure $\cJwh$ on 
$\Xwh = tot(V^{\vee})$ is integrable.

Recall that the choice of a connection $\nabla$ gives rise to a splitting 
$(D,\al)$:

\[
\xymatrix@1{0 \ar[r] & \ar[r] \pi^*V \ar[r]^-{j}  & T_X \ar@/^0.5pc/[l]^-{D} 
\ar[r]^-{d\pi} &
\ar[r]^-{} \ar@/^0.5pc/[l]^-{\al} \ar[r] \pi^*T_M \ar[r] & 0},
\]

\noindent
of the tangent sequence of $X \to M$.

Now if that {\em $\nabla$ is flat} it is known \cite{Kobayashi} that
we may find in a neighborhood of any point of $M$ a frame, $\{ e_i \}$ 
such that $\nabla e_i = 0$.  We will call $\{ e_i \}$ a {\em flat frame}.  
Given a flat frame, along with the corresponding vertical coordinates 
$\{\xi_i\}$ we have
that for any choice of coordinates $x_i$ on the base, the functions $\xi_i$ 
together with 
$y_i = x_i \circ \pi$ form a coordinate system on $X$ and  
$\alpha( \pi^* \pdxi) = (1 - j \circ D) \pdyi = \pdyi$ follows from 
the expression in this frame (see \ref{eqn:Ddef}) for $D$.  We define 
a frames $\{ f_i \}$ for $\pi^*V$ and $\{f^i\}$ for $\pi^*V^{\vee}$ by 
using the pullbacks $f_i = \pi^*e_i$ and $f^i = \pi^*e^i$, where $\{e^i\}$ 
is a dual frame to $\{e_i\}$.  

\begin{rem}
Notice that $\nabla$ is flat if and only if the image of 
$\alpha$ is involute.  Hence in this case we have 
a horizontal foliation instead of just a horizontal 
distribution.  We will consider the geometry of a pair of 
transverse foliations and its interaction with a generalized complex structure 
in section \ref{s:transverse}.
\end{rem}

\begin{defin}\label{def:semiflat}
If $\nabla$ is a {\em flat} connection on a rank $n$ vector bundle $V$ over a 
real $n$-manifold then a {\em $\nabla$-semi-flat} generalized almost 
complex structure on 
$X = tot(V)$ is an adapted, $\nabla$-lifted 
(see \ref{def:nablaadapted}) generalized almost complex structure.
\end{defin}

Let $\cS$ be the sub-sheaf of flat sections of $V \oplus V^{\vee}$.  Consider the isomorphism of vector bundles 

\begin{equation}\label{e:Mequation1}
\cM :V \oplus V^{\vee}  \to T_M \oplus T_M^{\vee},  \ \ \ \
\cM=\left(
\begin{array}{cc}
\cJ_{13}         &   \cJ_{22}^{\vee} \\
-\cJ_{31}^{\vee} &   \cJ_{12}^{\vee} \\
\end{array}
\right),  \ \ \ \
\cM^{-1}=\left(
\begin{array}{cc}
-\cJ_{12}         &   -\cJ_{22} \\
-\cJ_{31}          &   \cJ_{13}^{\vee} \\
\end{array}
\right).
\end{equation}

\noindent
With this notation we have

\begin{thm} \label{thm:integrability}
If $V$ is a vector bundle on a manifold $M$, then a semi-flat generalized
almost complex structure $\cJ = F^{-1} (\pi^*\cK) F$ on $X = tot(V)$ is 
integrable if and only if all pairwise Courant brackets of sections 
of the sheaf $\cM (\cS)$ vanish.
\end{thm}

Notice that this condition is expressed entirely in terms of data 
on the base manifold $M$.  Furthermore, we will see that this 
theorem implies the following 
corollary.

\begin{cor} \label{cor:mirror}
The generalized almost complex structure  $\cJ = F^{-1} (\pi^*\cK) F$ 
on $X = tot(V)$ is integrable if and only if the generalized almost
complex structure  $\cJwh = \Fiwh (\pih^*\cKwh) \Fwh$ on $\Xwh = tot(V^{\vee})$
is integrable.  
\end{cor}

\begin{rem}
In other words the mirror symmetry transformation is 
a bijective correspondence between 
$\nabla$-semi-flat generalized complex structures on $X$ and 
$\nabla^{\vee}$-semi-flat 
generalized complex structures on $\Xwh$.
\end{rem}

\begin{example}
If $\nabla$ is any flat, torsion-free connection on $T_M$, we can put 
a canonical complex structure on $tot(T_M)$.  See section \ref{s:examples} 
for more details.  This construction was first done in \cite{Dom}.  It is easy 
to see that the mirror structure is {\em always} the canonical symplectic 
structure on $tot(T_M^{\vee})$.
\end{example}

\noindent
{\bf Proof of Theorem \ref{thm:integrability}.}

Let us analyze the condition that the 
$(+\smo)$ eigenbundle $E$ be involute. The bundle $E$ is the graph 
of the isomorphism

\[
-\smo \cJ\big\vert_{image(j \oplus D^{\vee})\otimes \mathbb{C}}
: {image(j \oplus D^{\vee})\otimes \mathbb{C}} \to  {image(\al \oplus d\pi^{\vee})\otimes \mathbb{C}}.
\]

\noindent  
It suffices to analyze involutivity it locally on 
the base manifold.  Note that in the local 
frame and coordinates which we have chosen, we have the following formulas.

\

\ \ \ \ \ \ \ \ \ $j(f_i) = \pdxii$\hspace{25 mm}$D(\pdxii) = f_i$\hspace{25 mm}$D(\pdyi) = 0$

\

\ \ \ \ \ \ \ \ \ $\al(\pdxi) = \pdyi$\hspace{17 mm}$d\pi(\pdyi) = \pi^*\pdxi$ \hspace{10 mm}$d\pi(\pdxii) = 0$

\

\ \ \ \ \ \ \ \ \ $D^{\vee}(f^i) = d\xi_i$\hspace{25 mm}  
$\al^{\vee}(dy_i)= \pi^*dx_i$\hspace{20 mm}$\al^{\vee}(d\xi_i)= 0$

\

\ \ \ \ \ \ \ \ \ $(d\pi)^{\vee}(\pi^{*}dx_i)= dy_i$\hspace{13 mm}
$j^{\vee}(d\xi_i) = f^i$\hspace{27 mm}
$j^{\vee}(dy_i) = 0$

\

Furthermore, an isotropic sub-bundle of 
$(T_X \oplus T_X^{\vee})\otimes \mathbb{C}$ is involute 
if and only if it has a basis of sections whose pairwise Courant brackets 
are themselves sections of the original bundle.  This follows immediately
from the Leibniz property of the Courant Bracket see e.g. \cite{Xu, Cou}.  
This property says that 
\begin{equation} \label{e:intfrombasis}
\bigl[v_1 +\al_1, f(v_2 + \al_2) \bigr] = 
f \bigl[v_1 +\al_1,v_2 + \al_2 \bigr] 
+ v_1(f)(v_2 + \al_2) + \pair{v_1 +\al_1}{v_2 + \al_2} df
\end{equation}
for all vector fields $v_1$ and $v_2$, one-forms $\al_1$ and $\al_2$ 
and functions $f$. 
Let $U$ is the coordinate neighborhood of the base.  
We will analyze involutivity in $\pi^{-1}(U)$, using the coordinate system 
and frame described above.  Involutivity of $E$ is equivalent to 
the condition that 
$\bigl[a_i, a_j \bigr]$, $\bigl[a_i, b_j \bigr]$, 
and $\bigl[b_i, b_j \bigr]$ are all sections of $E$ where 
\[a_i = j(f_i) - \smo \cJ(j(f_i))
\] \noindent
and \[b_i = D^{\vee}(f^i) -\smo \cJ(D^{\vee}(f^i)).
\] \noindent
Using the special form of $\cJ$ we have 
\[a_i = \pdxii -\smo \al \pi^{*} (\cJ_{13} e_i) 
+ \smo (d\pi)^{\vee} \pi^{*} (\cJ_{31}^{\vee} e_i)
\] \noindent 
and 
\[b_i =  d \xi_i +\smo \al \pi^{*}(\cJ_{22}^{\vee} e^i)
+\smo d \pi^{\vee} \pi^{*}(\cJ_{12}^{\vee} e^i)
\]
\noindent
Hence we have that  
\settowidth{\dummy}{$\bigl[ a_i, a_j \bigr] \;$} 
\begin{align*}
\bigl[ a_i, a_j \bigr] 
& =  \bigl[ \pdxii -\smo \al \pi^{*} (\cJ_{13} e_i) 
+\smo (d\pi)^{\vee} \pi^{*} (\cJ_{31}^{\vee} e_i), 
\pdxij -\smo \al \pi^{*} (\cJ_{13} e_j)  
+\smo (d\pi)^{\vee} \pi^{*} (\cJ_{31}^{\vee} e_j)\bigr] 
\end{align*}
\begin{align*}
\hspace{\dummy} & =  [\pdxii -\smo \al \pi^{*} (\cJ_{13} e_i),
  \pdxij -\smo \al \pi^{*} (\cJ_{13} e_j)] \\ 
& \quad + \myiota_{\pdxii -\smo \al \pi^{*} (\cJ_{13} e_i)} 
d \smo (d\pi)^{\vee} \pi^{*} (\cJ_{31}^{\vee} e_j) 
- \myiota_{\pdxij -\smo \al \pi^{*} (\cJ_{13} e_j)} 
d \smo (d\pi)^{\vee} \pi^{*} (\cJ_{31}^{\vee} e_i) \\
& \quad + (1/2) d \myiota_{\pdxii -\smo \al \pi^{*} (\cJ_{13} e_i)}  
\smo (d\pi)^{\vee} \pi^{*} (\cJ_{31}^{\vee} e_j)
- (1/2) d \myiota_{\pdxij 
-\smo \al \pi^{*} (\cJ_{13} e_j)} \smo (d\pi)^{\vee} \pi^{*}
(\cJ_{31}^{\vee} e_i) \\
& = -[\al \pi^{*} (\cJ_{13} e_i), \al \pi^{*}
  (\cJ_{13} e_j)] \\ 
& \quad + \myiota_{\al \pi^{*} (\cJ_{13} e_i)} 
d {(d\pi)}^{\vee} \pi^{*} (\cJ_{31}^{\vee} e_j) 
- \myiota_{\al \pi^{*} (\cJ_{13} e_j)} 
d {(d\pi)}^{\vee} \pi^{*} (\cJ_{31}^{\vee} e_i) \\
& \quad + (1/2) d \myiota_{\al \pi^{*} (\cJ_{13} e_i)} 
(d\pi)^{\vee} \pi^{*} (\cJ_{31}^{\vee} e_j)
- (1/2) d \myiota_{ 
\al \pi^{*} (\cJ_{13} e_j)} (d\pi)^{\vee} \pi^{*} (\cJ_{31}^{\vee}
e_i) \\
& = -\al[\pi^{*} (\cJ_{13} e_i), \pi^{*} (\cJ_{13} e_j)] \\
& \quad + \myiota_{\al \pi^{*} (\cJ_{13} e_i)} 
 {(d\pi)}^{\vee} \pi^{*} d(\cJ_{31}^{\vee} e_j) 
- \myiota_{\al \pi^{*} (\cJ_{13} e_j)} 
 {(d\pi)}^{\vee} \pi^{*} d(\cJ_{31}^{\vee} e_i) \\
& \quad + (1/2) d \pi^{*}\myiota_{\cJ_{13} e_i} (\cJ_{31}^{\vee} e_j) 
- (1/2) d \pi^{*}\myiota_{\cJ_{13} e_j} (\cJ_{31}^{\vee} e_i) 
\end{align*}
\begin{align*}
\hspace{\dummy} & = -\al\pi^{*}[\cJ_{13} e_i, \cJ_{13} e_j] \\
& \quad + {(d\pi)}^{\vee} \pi^{*} \myiota_{\cJ_{13} e_i} 
d(\cJ_{31}^{\vee} e_j) 
- {(d\pi)}^{\vee} \pi^{*} \myiota_{\cJ_{13} e_j} 
 d(\cJ_{31}^{\vee} e_i) \\
& \quad + (1/2)(d\pi)^{\vee} \pi^{*} d \myiota_{\cJ_{13} e_i} 
(\cJ_{31}^{\vee} e_j)  
- (1/2) (d\pi)^{\vee}  \pi^{*} d \myiota_{\cJ_{13} e_j}
(\cJ_{31}^{\vee} e_i) \\
& =  - (\al +(d\pi)^{\vee})\bigl[(\pi^*\cJ_{13})e_i -
  (\pi^*\cJ_{31}^{\vee})e_i,             \pi^{*}(\cJ_{13} e_j) -
  \pi^{*}(\cJ_{31}^{\vee})e_j) \bigr] \\ 
& =  - (\al +(d\pi)^{\vee}) \pi^{*}\bigl[\cJ_{13}e_i -
  \cJ_{31}^{\vee}e_i, \cJ_{13} e_j  
- \cJ_{31}^{\vee}e_j \bigr] 
\end{align*}
\noindent
or
\begin{equation}\label{e:abrack}
\bigl[ a_i, a_j \bigr] = - (\al +(d\pi)^{\vee}) 
\pi^{*}\bigl[\cJ_{13}e_i - 
\cJ_{31}^{\vee}e_i, \cJ_{13} e_j - \cJ_{31}^{\vee}e_j \bigr]
\end{equation}

\noindent
Similarly we have 

\begin{equation}\label{e:abbrack}
\bigl[a_i, b_j\bigr] =  (\al +(d\pi)^{\vee}) \pi^* \bigl[\cJ_{13} e_i
  - \cJ_{31}^{\vee}e_i , \cJ_{22}^{\vee} e^j + \cJ_{12}^{\vee} e^j
  \bigr] 
\end{equation}

\noindent
and 

\begin{equation}\label{e:bbbrack}
\bigl[b_i, b_j\bigr] = -(\al +(d\pi)^{\vee}) \pi^*
\bigl[\cJ_{22}^{\vee} e^i + \cJ_{12}^{\vee} e^i , \cJ_{22}^{\vee} e^j
  + \cJ_{12}^{\vee} e^j \bigr]. 
\end{equation}

\

\noindent
The right hand sides of all three of these expressions are sections 
of the vector bundle $image(\al + (d\pi)^{\vee})\otimes{\mathbb{C}}$.  
Therefore, the right hand sides are sections of $E$ and 
in particular be sections of the graph of a map of vector bundles from 
$image(j +D^{\vee}) \otimes \mathbb{C}$ to \linebreak
$image(\al +(d\pi)^{\vee})\otimes \mathbb{C}$ if and only if 
$\bigl[a_i, a_j\bigr]$, 
$\bigl[a_i, a_j\bigr]$, and $\bigl[b_i, b_j\bigr]$ all vanish.  
This is precisely the statement of \ref{thm:integrability}: that 
all pairwise Courant brackets between sections 
of $\cM(\cS)$ vanish. 

\ \hfill $\Box$

\bigskip

Notice now that if we replace the vector bundle $V$ by $V^{\vee}$ 
and $\cK$ by $\cKwh$ (see \ref{e:K_check}, \ref{rem:Jswitch}) then 
$\cM$ gets replaced by  

\begin{equation}\label{e:Mequation2}
\widehat{\cM}=\left(
\begin{array}{cc}
-\cJ_{22}^{\vee}         &   -\cJ_{13}      \\
-\cJ_{12}^{\vee}         &   \cJ_{31}^{\vee} \\
\end{array}
\right).
\end{equation}

\noindent
but $\cM (\cS) = \widehat{\cM} (\cS^{\vee})$.  Therefore
we have also proven \ref{cor:mirror}.  \hfill $\Box$

\

\bigskip

It is also clear from this proof and using equation (\ref{e:intfrombasis}), 
that if $\cJ$ is integrable, then the two almost Dirac
structures $\Delta = \cK(V) = \cKwh(V)$ and 
$\widehat{\Delta} = \cKwh(V^{\vee}) = \cK(V^{\vee})$ are as well.  The 
vector bundle $\Delta$ inherits the same flat connection from $V$ via 
$\cK$ or $\cKwh$.  Similarly, $\widehat{\Delta}$ inherits the same 
flat connection from $\cK$ or $\cKwh$.

\begin{cor}\label{cor:Dirint}
For a flat connection $\nabla$ on a vector bundle $V$ over $M$, a 
$\nabla$-semi-flat generalized complex structure on the total space 
of $V$ induces a pair of transverse Dirac structures on $M$.  These Dirac 
structures inherit flat connections.
\end{cor}   \hfill $\Box$

\

\begin{rem}
The geometry of a pair of transversal Dirac subbundles was recently studied 
by A. Wade and found to be equivalent to a generalized paracomplex structure 
as defined in \cite{Wade}.  Furthermore using the analysis of the  
integrability condition in terms of local systems above, the two 
Dirac structures that we have identified above form a {\em pair of Dirac 
structures} (see e.g. \cite{Dorfman}) in the sense of Gelfand and Dorfman and 
therefore leads to method of constucting integrable hierarchies with respect 
to the two Poisson structures comming from the two Dirac structures.  This 
remark also applies to the torus bundle case below.  Another 
overlap with the mathematics of integrable systems is also noted in 
section \ref{s:transverse} and these overlaps will be the subject of 
future work.    
\end{rem}

Note that a {\em generalized K\"{a}hler structure} is defined \cite{Gua} to be 
a generalized almost K\"{a}hler structures where both of the two 
generalized almost complex structures are integrable.  Therefore 
we have also proven (using the results of subsection \ref{ss:genkah}) 
the following.

\begin{cor}\label{cor:Kint}
The correspondence in Corollary \ref{cor:mirror}, gives  
a bijective correspondence between $\nabla$-semi-flat 
generalized K\"{a}hler structures on $X$ and $\nabla^{\vee}$-semi-flat 
generalized K\"{a}hler structures on $\Xwh$.
\end{cor}   \hfill $\Box$

\section{From vector bundles to torus bundles}\label{s:tbundle}

In this section we describe generalized complex structures on (real) 
torus bundles 
with sections and their mirrors.  The base of our torus 
bundles will turn out to 
support a pair of complimentary Dirac structures.  

\subsection{The geometry of torus bundles and the dual of a torus
  bundle}\label{ss:torbungeom}  

Let $\cX \rar{p} M$ be a fiber bundle over a real manifold $M$ for 
which the fibers have the diffeomorphism type of a real torus of 
dimension $n$.  We call this a {\em torus bundle}.  So we have 
for $U$ small in the base, local 
isomorphisms of fiber bundles $p^{-1}(U) \cong U \times \boldsymbol{T}$.  
Assume that this fiber bundle possesses a global (smooth)
section $s$.  This is equivalent to assuming that the structure group of 
the bundle is $Diff(\boldsymbol{T}, 0)$ as opposed to $Diff(\boldsymbol{T})$.  
However, since the connected component of $Diff(\boldsymbol{T}, 0)$ is 
contractible, the structure group of the bundle may be 
reduced to those diffeomorphisms which respect the group structure: 
$Aut(\boldsymbol{T}) \cong GL(n, \mathbb{Z})$.  Recently this issue was 
discussed in \cite{Jardim}.  We consider this to have 
been done and regard $s$ as the zero section.    We therefore consider $X$ 
as a (Lie) group bundle or bundle of (Lie) groups.  Recall that for a 
bundle of Lie groups modeled on a Lie group $\boldsymbol{G}$ 
(we sometimes call this simply a $\boldsymbol{G}$-bundle) 
we have local 
maps $\rho^{-1}(U) \cong U \times \boldsymbol{G}$ and the transition maps 
$U \cap V \times \boldsymbol{G} \cong U \cup V \times \boldsymbol{G}$ 
restrict to a Lie group isomorphism of $\boldsymbol{G}$ on each fiber.  

Consider the tangent sequence

\[
0 \rar{} T_{\cX/M} \rar{} T_{\cX} \rar{dp}
p^*T_M \rar{} 0
\]

As in case of vector bundles,  
$T_{\cX/M}$ is a pullback.  Indeed, let $V = s^*T_{\cX/M}$, then 
we have that $T_{\cX/M} \cong p^* V$.  This 
follows from the following simple observation.  

\begin{lem} \label{lem:reltan}
Let $\boldsymbol{G}$ be a Lie group and $\cY \rar{\rho} N$ a 
$\boldsymbol{G}$-bundle.  Call the zero-section $s$.  Then we have  
$\rho ^* s^* T_{\cY/N} \cong T_{\cY/N}$.
\end{lem}

\noindent
{\bf Proof.}

Write any section $\sigma$ of $(s \circ \rho)^{-1} T_{\cY/N}$ over 
$U \subseteq \cY$ as $\sigma = \sigma_0 \circ s \circ \rho$ where 
$\sigma_0$ is a section of $T_{\cY/N}$ over $s (\rho(U))$.  Now using 
the local group structure we may push $\sigma_0$ 
forward along the fibers.  
The transition maps respect the 
group structure of $\boldsymbol{G}$ and therefore these vector fields patch to a section 
of $T_{\cY/N}$ over 
$\rho^{-1} (\rho(U))$ and then we can restrict this section to $U$.  This 
gives a morphism of vector bundles 

\[
(s \circ \rho)^* T_{\cY/N} \rar{\psi} T_{\cY/N}.
\] 

Over a point $y \in \cY$, when 
we look in one of the trivial neighborhoods, 
$\rho^{-1}(U) \cong U \times \boldsymbol{G}$ where $y = (u, g)$, the map 
becomes just the obvious map $Lie(\boldsymbol{G}) \to T_g {\boldsymbol{G}}$ which is clearly 
an isomorphism.  Hence we can conclude that the map $\psi$ gives an 
isomorphism 
$\rho ^* s^* T_{\cY/N} \cong T_{\cY/N}$.    \hfill $\Box$

\

Furthermore, if we use the fact that the torus {\em compact and connected}, 
the sheaf of sections, $V$ of $V = s^*T_{\cX/M}$ is 
isomorphic to $(R^1 p_* \mathbb{R})^{\vee} \otimes C^\infty_M$.  

\begin{lem} \label{lem:Vflat}
If $\cX \rar{p} M$ is a $\boldsymbol{T} = Lie(\boldsymbol{T})/ \Gamma$ 
bundle with structure group $Aut(\Gamma)$ then 
$V \cong (R^1 p_* \mathbb{R})^{\vee} \otimes C^\infty_M$.
\end{lem}

\noindent
Notice that this is just a relative version of the natural isomorphism 
$Lie(\boldsymbol{T})^{\vee} \cong H^1(\boldsymbol{T}, \mathbb{R})$ which 
is described for example in \cite{Bred}.

\noindent
{\bf Proof.}

Notice that $V = s^* T_{\cX/M}$ is a $Lie(\boldsymbol{T})$-bundle on $M$ 
with structure group $Aut(\Gamma)$.  Let $\Lambda \subseteq tot(V)$ be the 
lattice induced by $\Gamma$ and $\cS_{\Lambda}$ be its sheaf of sections.  
There is a morphism of presheaves of abelian groups 
\[
\cS_{\Lambda} \to [U \mapsto H_1(p^{-1}(U), \mathbb{Z})].
\]  It is given (for $U$ connected) 
by sending 
$\lambda \in \cS_\Lambda (U)$ to the homology class of the image in $\cX$ of 
the line in $tot(V)$ connecting the point $0$ over $m$ to the point 
$\lambda(m)$ for some $m \in U$.  For $U$ small enough this is an 
isomorphism using the K\"{u}nneth theorem
\[
\cS_\Lambda (U) \cong H_1(p^{-1}(U), \mathbb{Z}) = 
Hom_{\mathbb{Z}}(H^1(p^{-1}(U), \mathbb{Z}), \mathbb{Z}).
\]  
Hence we have isomorphisms of sheaves 
\[
\cS_{\Lambda} \cong 
Hom_{\mathbb{Z}}(R^1 p_* \mathbb{Z}, \mathbb{Z})\mbox{,\   \ } 
\cS_{\Lambda} \otimes \mathbb{R} \cong (R^1 p_* \mathbb{R})^{\vee}, 
\mbox{  and \   \ }
V \cong (R^1 p_* \mathbb{R})^{\vee} \otimes C^\infty_M.
\]    \hfill $\Box$

Therefore we have a flat Gauss-Manin connection $\nabla$ on $V$ given as the 
image of $1 \otimes d$ under this isomorphism.  Recall that in the case 
of vector bundles we had for each flat connection on $V$ a 
(potentially) inequivalent mirror 
symmetry transformation.  By contrast, in the case of torus bundles, 
the topology of a torus bundle has given us a natural flat connection on $V$ 
and so we need not make any additional choices.

The multisection of $V$ given by $\Lambda = tot(\cS_\Lambda)$ acts on 
$X = tot(V)$ and
the orbits are the fibers of the natural map from $X$ to $\cX$.  Hence we have
a diffeomorphism $X/ \Lambda \cong \cX$.  Under this quotient, the 
multisection goes 
to $s(M)$.  The sheaf of sections of $\cX$ becomes the sheaf of groups 
$V/\cS_\Lambda 
\cong ((R^1 p_* \mathbb{R})^{\vee}) \otimes C^\infty_M)  
/ (R^1 p_* \mathbb{Z})^{\vee}$ where the
zero section has image $s(M)$.  The isomorphism 
$X/\Lambda \cong \cX$ is an isomorphism of $GL(n,\mathbb{Z})$ fiber bundles 
over $M$.

\

Now the 
{\em dual torus bundle} is defined to be 
$\cXwh = \Xwh / \widehat{\Lambda} \rar{\ph} M$ where
$\widehat{\Lambda} = tot({\cS_{\Lambda^{\vee}}})$ 
and ${\cS_{\Lambda^{\vee}}}= Hom_\mathbb{Z} (\cS_\Lambda, \mathbb{Z}) \subseteq
V^{\vee}$.  Furthermore, ${\cS_{\Lambda^{\vee}}} \cong R^1 \ph_* \mathbb{Z}$ 
and 
$V^{\vee} \cong R^1 \ph_* \mathbb{R} \otimes C^\infty_M$.  This gives
a flat connection on $V^{\vee}$ which is of course just the dual 
connection $\nabla^{\vee}$.  Also 
${\cS_{\Lambda^{\vee}}} \otimes_\mathbb{Z} \mathbb{R}$ is the sheaf of
flat sections of $V^{\vee}$ with respect to 
$\nabla^{\vee}$.  The sheaf of sections of $\cXwh$ over $M$ is a sheaf
of groups given by
$V^{\vee} / {\cS_{\Lambda^{\vee}}} \cong  
((R^1 \ph_* \mathbb{R})) \otimes C^\infty_M)  / (R^1 \ph_* \mathbb{Z})$.  
We then have a global section $\sh$ of $\cXwh$ over $M$ which is the 
zero section
and satisfies that $\sh(M)$ is the image of the multisection 
$\widehat{\Lambda}$ under
the quotient map.

We saw in section \ref{s:vectmirror} that if we let $X = tot(V)$,  
then $\nabla$ gives us a splitting 
$D$ of the tangent sequence of the map $X \rar{\pi} M$.  
We can use this to split the tangent sequence of the map $\cX \rar{p} M$.  
Consider the following diagram where we have decomposed $\pi$ as $p \circ q$.

\[
\xymatrix@1{X \ar[r]^-{q} & \cX
\ar[r]^-{p} & M \ar@/^0.5pc/[l]^-{s}},
\]
\noindent
We may push forward the exact sequence

\[
0 \rar{} \pi^*V \rar{} T_{X} \rar{d\pi} {\pi}^*T_M \rar{} 0
\]
\noindent
which is split by $\pi^*(V) \lar{D} T_{X}$ to the exact sequence

\[
0 \rar{} q_*\pi^*V \rar{} q_*T_{X} \rar{q_*{d\pi}} q_*{\pi}^*T_M \rar{} 0
\]
\noindent
which is split by $q_*\pi^*(V) \lar{q_*D} q_*T_{X}$.  Furthermore, 
$\Lambda$ naturally acts on all three of these sheaves and if we take 
the $\Lambda$ invariants of each term of this sequence we recover 
precisely the exact sequence that we want to split, namely:

\[
0 \rar{} p^*V \rar{} T_{\cX} \rar{dp} p^*T_M \rar{} 0
\]

Therefore, the only thing to check in order to split this sequence is that 
the map $D$ satisfies 
$(dt_\lambda)(Dw) = D((dt_\lambda)(w))$ where for some 
small $U \subseteq M$ and small $U' \subseteq p^{-1}(U)$
$w$ is a section of $T_X$ over $q^{-1}(U')$, $\lambda$ is a component 
of $\Lambda \cap \pi^{-1}(U)$ and $t_\lambda : X \to X$ is the action of 
addition of $\lambda$.  However, 
\[
(dt_\lambda) (Dw) = (dt_\lambda) (((\pi^*\nabla)S)w) 
= ((\pi^*\nabla)(S+\lambda))((dt_\lambda)w) 
= ((\pi^*\nabla)S)((dt_\lambda)w) = D((dt_\lambda)w)
\]
\noindent
due to the fact that that the sections of the lattice are flat.

\subsection{Generalized complex structures on torus bundles and the
  mirror transformation}\label{ss:gcplxtor} 

We will now use the same names as in the vector bundles case 
for the splittings of the tangent sequences of $\cX$ and $\cXwh$.  
That is:

\[
\xymatrix@1{0 \ar[r] & \ar[r] p^*V \ar[r]^-{j}  & T_{\cX} \ar@/^0.5pc/[l]^-{D} \ar[r]^-{dp} &
\ar[r]^-{} \ar@/^0.5pc/[l]^-{\al} \ar[r] p^*T_M \ar[r] & 0},
\]
\noindent
and

\[
\xymatrix@1{0 \ar[r] & \ar[r] \ph^*V^{\vee} \ar[r]^-{\jh}  & T_{\cXwh} \ar@/^0.5pc/[l]^-{\Dwh} 
\ar[r]^-{d\ph} &
\ar[r]^-{} \ar@/^0.5pc/[l]^-{\alh} \ar[r] \ph^*T_M \ar[r] & 0},
\]

\noindent
Since we will be using only one connection in the case of 
torus bundles, we will 
drop $\nabla$ from the notation.

\begin{defin} \label{def:semiflattor}
If $M$ is an $n-$dimensional real manifold and $\cX \to M$ is a real torus 
bundle with fiber dimension $n$ and zero section $s$ then we call 
a generalized almost complex 
structure $\cJ$ on $\cX$ which comes from (see section \ref{s:vectmirror})
 an adapted generalized almost 
complex structure $\cK$ on $s^*T_{\cX/M} \oplus T_M = V \oplus T_M$ 
a {\em semi-flat} generalized almost complex structure.
\end{defin}
\noindent
Recall that ``adapted'' just means that 
\[
\cK(s^*T_{\cX/M} \oplus s^*T_{\cX/M}^{\vee}) = T_M \oplus T_M^{\vee}
\]

As in the vector bundle case, there is a bijective correspondence 
between semi-flat generalized almost complex structures on $\cX$ and
$\cXwh$.  The proof is precisely the same, except the choice of local 
coordinates is local along the base and the fiber, instead of just along 
the base.

\begin{thm} \label{thm:torintegrability}
A semi-flat generalized almost complex structure $\cJ$ on a torus bundle 
$\cX \to M$ with zero section $s$ is integrable if and only if 
\[
\bigl[\cK (\cS \oplus \cS^{\vee}),\cK (\cS \oplus \cS^{\vee}) \bigr] = 0
\]  
where $\cS$ is the sheaf of flat sections of $s^*T_{\cX/M}$.
\end{thm}       \hfill $\Box$ 

\begin{cor} \label{cor:tormirror}
A semi-flat generalized almost complex structure $\cJ$ on a torus bundle 
${\cX \to M}$ is integrable if and only if its mirror structure $\cJwh$ on the 
dual torus bundle $\cXwh \to M$ is integrable.
\end{cor}       \hfill $\Box$ 

\begin{rem}
This means that we have given a bijective correspondence 
between semi-flat generalized complex structures on 
$\cX$ and semi-flat generalized complex structures on $\cXwh$.  The 
same holds true for generalized K\"{a}hler structures in which 
both of the generalized complex structures are semi-flat.
\end{rem}

\begin{cor} \label{cor:tordirac}
A semi-flat generalized almost complex structure $\cJ$ on a torus bundle 
${\cX \to M}$ induces a pair of almost Dirac structures 
\[
\Delta = \cK(s^*T_{\cX/M}),  \  \  
\widehat{\Delta} = \cK(s^*T_{\cX/M}^{\vee}) \subseteq T_M \oplus T_M^{\vee}
\]
Each carries its own flat connection and these Dirac structures are exchanged 
under mirror symmetry.  If $\cJ$ is integrable then $\Delta$ and 
$\widehat{\Delta}$ are integrable.  
\end{cor}      \hfill $\Box$

Now in the case when the generalized complex structure on the torus bundle 
is of symplectic type and the torus fibers are Lagrangian this result 
reproduces the starting point of the work \cite{Pol} where the torus bundle 
is written as $tot(T_M^{\vee})/\Lambda$ and $\Delta$ is the 
Dirac structure $T_M^{\vee}$, which inherits a flat connection $\nabla$.  
The mirror manifold 
$tot(T_M)/\Lambda^{\vee}$ inherits a complex structure as 
explained in \cite{Pol} constructed using the dual connection 
$\nabla^{\vee}$ which is both flat and torsion-free.  This corresponds 
to the canonical almost complex structure on $tot(T_M)$ associated 
to a connection on $T_M$ which is known \cite{Dom} to be integrable 
if and only if the connection is flat and torsion-free.    

\section{The cohomology of torus bundles}\label{s:torcoho}
  
Consider the diagram

\[
\xymatrix{
& \cX\times_{M} \cXwh \ar[dl]_-{\qh} \ar[dr]^-{q}
 & \\
\cX \ar[dr]_-{p} & & \cXwh
\ar[dl]^-{\ph} \\
& M &
}
\]

\noindent
Now the space $\cX \times_{M} \cXwh$ is endowed with a global closed 
two form given as 
$\Xi = \frac{1}{2 \pi \smo} \mathcal{F}$ where $\mathcal{F}$ 
is the curvature of 
a connection on the relative 
Poincar\'{e} (line) bundle.  See \cite{Pioli} for an explanation of the 
relative Poincar\'{e} bundle in this context.  Now we would like to introduce 
a relative version of a map given \cite{Orlov} in the context of 
mirror symmetry of 
abelian varieties as introduced by Mukai \cite{Muk}.  This idea 
has appeared in various ways in \cite{Allessandro, Van, Gua, Mathai} and 
the references therin.

\begin{lem} \label{lem:cohomology}
If the bundle $\qh^{*} T_{\cX/M}$ is orientable  
then we have a morphism (independent of the choice of orientation) 
of sheaves of $C^\infty_M$ modules
$p_* \Omega^{\dot}_{\cX} \to \ph_* \Omega^{\dot}_{\cXwh}$ is a 
morphism of the de Rham complexes.  Therefore 
this morphism gives a map of presheaves 
\[
[U \mapsto H^{\dot}(p^{-1} (U),\mathbb{R})] 
\to [U \mapsto H^{\dot}(\ph^{-1} (U),\mathbb{R})].
\]
This map of presheaves induces an 
isomorphism of the sheafifications 
$R^{\dot}p_{*} \mathbb{R} \to  R^{\dot}\ph_{*} \mathbb{R}$ which decomposes
into isomorphisms $R^{j}p_{*} \mathbb{R} \to  R^{n-j}\ph_{*} \mathbb{R}$ 
for $j=1,\ldots,n$.
\end{lem}

\noindent
{\bf Proof.}

We have a map $\qh^{\flat} = (d\qh)^{\vee} \circ \qh^*$ from 
$p_*\Omega^{j}_{\cX}$ to 
$p_* \qh_* \Omega^{j}_{\cX \times_M \cXwh}$  given by pulling back
differential forms.  Clearly, $\qh^{\flat}$ commutes with 
the de Rham differentials.  
Observe that the map $q$ makes $\cX \times_M \cXwh$ into a torus bundle 
over $\cXwh$.  The relative tangent bundle of the tangent sequence of the 
map $q$ is isomorphic to $\qh^*T_{\cX/M}$.  Therefore we also have a 
map $q_*$ which integrates along 
the fibers and maps 
$q_* \Omega^{k}_{\cX \times_M \cXwh}$ to $\Omega^{k-n}_{\cXwh}$.  
Explicitly, 
if we take our global section $s$ over $\cX \times_M \cXwh$ of 
$\Wedge^n \qh^*T_{\cX/M}$ and the 
corresponding global section $t$ of  $\Wedge^n \qh^*T_{\cX/M}^{\vee}$ then
$q_*(\Upsilon) = \int_{(\cX \times_M \cXwh) / \cX} ((\myiota_s \Upsilon) \wedge t)$.  This does not depend 
on the choice of $s$ but we do need the fibers of $q$ to be orientable 
manifolds to integrate over them.  
Since the
torus fibers of $\qh$ are manifolds without boundary we have that 
$q_*$ and also
its pushforward, $\ph_* [q_*]$ commutes with the de Rham differentials.
Now we can define the map 
$\mathbf{F.T.} : p_* \Omega^{\dot}_{\cX} \to \ph_* \Omega^{\dot}_{\cXwh}$ by 

\[
\mathbf{F.T.}(\mu) = \ph_* [q_*] (\qh^{\flat}(\mu) \wedge \exp(\Xi))  
\]

Now since $d \Xi =0$ we have 
$\mathbf{F.T.}(d \mu) = d \mathbf{F.T.}(\mu)$ and hence we get a 
map of presheaves $[U \mapsto H^{\dot}(p^{-1} (U),\mathbb{R})] 
\to [U \mapsto H^{\dot}(\ph^{-1} (U),\mathbb{R})]$.  In particular we have
a natural $\mathbb{R}-$linear map 
$H^{\dot}(\cX, \mathbb{R}) \to H^{\dot}(\cXwh, \mathbb{R})$.  In order to 
write the map on differential forms locally on the base, chose a trivializing
open neighborhood $U \subseteq M$ and $\mu \in \Omega^c(p^{-1}(U))$.  Then let
$\xi_i$ be the flat vertical coordinates on $p^{-1}(U)$ and $\eta_i$ be
the dual flat vertical coordinates on $\ph^{-1}(U)$.  In these 
coordinates we may assume without loss of generality that  

\[
s = {\partial/{\partial \xi_n}} \wedge \cdots \wedge {\partial/{\partial \xi_1}}\]

\noindent
on $p^{-1}(U)$.  (Every section may be extended to a global section.)

\noindent
Let us now express $\mu$ in local coordinates.

\[
\mu = \sum_{|J|=1, \ldots, c} {f_{J}} \Theta_J 
\wedge d\xi_{j_1} \wedge \cdots \wedge d\xi_{j_b}
\]

\noindent
Here, the $\Theta_J$ are pullbacks of  $(c-b)-$ forms from the base, 
$J = (j_1 , \ldots , j_b)$ where $j_1 < \cdots < j_b $ and
$f_J$ are functions on $p^{-1}(U)$.  A simple calculation shows that

\[\hat{\mu} = \mathbf{F.T.}(\mu) = \int_{\boldsymbol{T}} \mu 
\wedge \exp(d\xi_i \wedge d\eta_i) 
\]

\noindent
is given by

\[
\hat{\mu} = 
\sum_{|J|=1, \ldots, c} (-1)^{k_1 + \cdots + k_{n-b}} \theta_J 
\wedge d\eta_{k_1} \wedge \cdots \wedge d\eta_{k_{n-b}}
\int_{\boldsymbol{T}}  f_{J} d\xi_1 \wedge \cdots \wedge d\xi_n
\]

\noindent
where $k_1 < \cdots < k_{n-b}$ is the compliment to $J$.  Now suppose 
that $\mu$ is closed and that we consider the cohomology class 
$[\mu] \in H^j(p^{-1}(U), \mathbb{R})$, Using the K\"{u}nneth isomorphism

\[
H^c (p^{-1}(U), \mathbb{R}) \cong 
\bigoplus_{b=0,\ldots,c} H^{c-b}(U,\mathbb{R}) 
\otimes H^{b}(\boldsymbol{T}, \mathbb{R}) \cong
H^{c}(\boldsymbol{T}, \mathbb{R})
\]

\noindent
we may absorb the $f_J$ into the $\Theta_J$ in the above expression
and therefore since the cohomology of the tori $\boldsymbol{T}$ and 
$\boldsymbol{T}^{\vee}$ are generated by the classes 
$[d\xi_1 \wedge \cdots \wedge d\xi_j]$ and $[d\eta_1 \wedge \cdots \wedge d\eta_k]$ 
respectively we conclude that the map $\mathbf{F.T.}$ induces 
isomorphisms $R^{j}p_{*} \mathbb{R} \to  R^{n-j}\ph_{*} \mathbb{R}$ 
for $j=1,\ldots,n$ as promised.     \hfill $\Box$ 

\begin{cor}\label{cor:torspinors}
If $\cJ$ is a semi-flat generalized almost complex structure on a $n$-torus 
bundle with section $\cX$ on an $n$-manifold $M$ with 
associated spinor line bundle 
$L \subseteq \Wedge^{\dot} T_{\cX}^{\vee}\otimes{\mathbb{C}}$, and $\Lwh$ is 
the line bundle associated to the mirror structure $\cJwh$ on $\cXwh$, 
then 
\[
\mathbf{F.T.}(p_*L) = \ph_*\Lwh \subseteq \ph_* \Wedge T_{\cXwh}^{\vee}\otimes{\mathbb{C}}.
\]
\end{cor}

\noindent
{\bf Proof.}

This follows from tensoring the previous 
lemma with $\mathbb{C}$ and using 
lemma \ref{lem:fouriermukai}.      \hfill $\Box$ 

\begin{rem} The Fourier-Mukai transformation for spinors, combined with 
the formulae given in Lemma \ref{lem:fouriermukailem} can easily be used 
to show again that integrability is preserved by the 
mirror transformation we have 
desccribed.  As we have already shown this, we do not 
demonstrate it again with spinors.
\end{rem}

\begin{example}
Let $M = S^1$ or $\mathbb{R}$, $Z = V/\Lambda \times M$, 
$\widehat{Z}= V^{\vee}/\Lambda^{\vee} \times M$, where $V/\Lambda \cong S^1$.  
Let $x, \theta, \hat{\theta}$ be ``coordinates'' on $M$, $V/\Lambda$ and 
$V^{\vee}/\Lambda^{\vee}$ respectively.  Then for $f$ a complex valued 
smooth nowhere vanishing function on $M$. 

\[
\mathbf{F.T.}(e^{f d\theta \wedge dx}) 
= \int_{V/\Lambda} e^{f d\theta \wedge dx + d{\theta} \wedge d\hat{\theta}} 
= \int_{V/\Lambda} d\theta \wedge(f dx + d\hat{\theta}) = d\hat{\theta}+ f dx
\]
\noindent
When we take $f = i$, we see the spinor corresponding to a symplectic 
structure going to one representing a complex structure.
\end{example}

\begin{rem} Let $Z$ an $n$-torus bundle over a compact 
connected $n$-manifold such that $\qh^{*} T_{\cX/M}$ is orientable.    
Consider the ``Moduli-Space'' $SFGCY(\cX)$ of 
{\em semi-flat generalized Calabi-Yau structures}.  These are semi-flat 
generalized 
complex structures which are {\em generalized Calabi-Yau} \cite{H2}, meaning 
that the associated spinor line bundles $L$ have nowhere 
vanishing, closed, global 
sections.  

Since these sections are known \cite{Che,Gua, H2} to be either even 
or odd we may consider a ``period map'' \cite{H3, Huy} from this space into 
$\mathbb{P}(H^{even}(\cX, \mathbb{C})) \coprod \mathbb{P}(H^{odd}(\cX, \mathbb{C}))$.  Note that we are assuming here that for a fixed structure, 
different closed, nowhere vanishing global sections of $L$ 
define the same cohomology class up to multiplication by constants.  
Under this assumption, we have shown the existence of a commutative diagram.

\[
\xymatrix{
SFGCY(\cX) \ar[r] \ar[d] &  \mathbb{P}(H^{even}(\cX, \mathbb{C})) \coprod \mathbb{P}(H^{odd}(\cX, \mathbb{C})) \ar[d] \\
SFGCY(\cXwh) \ar[r] &
\mathbb{P}(H^{even}(\cXwh, \mathbb{C})) \coprod \mathbb{P}(H^{odd}(\cXwh, \mathbb{C})) 
}
\]
\end{rem}

\noindent
In the case that the torus bundles are $Z$ and $\Zwh$ trivial, 
the vertical map takes 
horizontal $i-$th cohomology to itself  
and vertical $i$-th cohomology to vertical $(n-i)$-th cohomology 
(both with multiplication by signs).

\begin{conj}
Let $(\cX, \cJ)$ be a compact generalized Calabi-Yau manifold of (real) 
dimension $2n$.  As we 
have mentioned in the previous remark, it would be desirable to know that 
there is a unique element in $\mathbb{P}(H^{\dot}(\cX, \mathbb{C}))$ 
associated with $\cX$.  Without this knowledge, the previous diagram would 
have to be modified by the appropriate restrictions on the left hand side.  
Therefore we conjecture (without overwhelming evidence) that if 
$\phi$ is a global, 
closed, nowhere vanishing differential form, representing $\cJ$, and $f$ is 
a nowhere zero smooth complex valued function such that 
$d(f \phi) = 0$, that $f$ is constant. If we call generalized Calabi-Yau 
manifolds satisfying this condition {\em Liouville} then it is easy to see 
that all symplectic manifolds are Liouville (take $\phi = e^{-i\om}$), 
compact Calabi-Yau manifolds are 
Liouville (take $\phi$ to be a nowhere zero holomorphic $n-$form), 
and products and $B-$field transformations of Liouville generalized
Calabi-Yau manifolds are Liouville. 
\end{conj}

\section{Transverse foliations and generalized K\"{a}hler
  geometry}\label{s:transverse} 

In this section we study in the abstract some of essential 
geometric details of our construction without reference to the 
specific context (e.g. the type of bundle).

\begin{defin} \label{def:compatible}
Suppose that $X$ is a foliated manifold and let $\cP \subseteq T_{X}$ be
the involute sub-bundle tangent to the leaves of the foliation.  We say
that a generalized complex structure $\cJ$ and $\cP$ are {\em
compatible} if there exists a complementary sub-bundle $\cQ \subseteq
T_{X}$ so that

\[
\cJ_{\big\vert \cP\oplus \op{Ann}(\cQ)} \; : \cP\oplus \op{Ann}(\cQ) \to \cQ
\oplus \Ann(\cP)
\]
is an isomorphism of vector bundles.  Under this condition, we will call
$\cQ$ a $\cJ$-compliment to $\cP$.
For $\cQ$ a $\cJ$-compliment to $\cP$
we will often tacitly identify ${\cP}^{\vee}$ 
with $\Ann{\cQ}$, and ${\cQ}^{\vee}$
with $\Ann{\cP}$.

\
\end{defin}

Notice that the $(+\smo)$ eigenbundle, $E$ of $\cJ$ is in this case 
necessarily transverse 
to both ${\cP}_{\bC} \oplus \Ann(\cQ)_{\bC}$ and  
${\cQ}_{\bC} \oplus \Ann({\cP})_{\bC}$.  
Hence $E$ is the graph of a map from  ${\cP}_{\bC} \oplus \Ann({\cQ})_{\bC}$ 
to  ${\cQ}_{\bC} \oplus \Ann({\cP})_{\bC}$.
In fact, it is easy to 
see that we have $E = graph(-\smo\cJ)$ where we consider $(-\smo\cJ)$ as 
a map from
${\cP}_{\bC} \oplus \Ann({\cQ})_{\bC}$ to 
${\cQ}_{\bC} \oplus \Ann({\cP})_{\bC}$.

\begin{examples} \label{exs:lagrangian} {\bf (i)} Suppose that
  $X$ is a manifold equipped with an involute
  distribution $\cP \subseteq T_{X}$ of half the dimension of $X$. 
  Let $\cJ$ be the generalized
  almost complex structure on $X$ corresponding to a non-degenerate real 
  two form $\omega$. Then $\cP$ and $\cJ$ are compatible if and only if $\cP$
  defines a Lagrangian foliation on $X$.  Indeed, the compatibility shows
  that $\omega$ defines an isomorphism from $\cP$ to $\Ann{\cP}$, which shows 
  that $\cP$ is Lagrangian.  Conversely, if $\cP$ is Lagrangian, then 
  by (see \cite{dS}) choosing an almost complex 
  structure $J$ so that the isomorphism
  $-\omega J:T_X \to T_X^{\vee}$ represents a Riemannian metric on $X$,
  it is easy to see that the vector bundle $J \cP$ is a $\cJ$-compliment to
  $\cP$, and so $J \cP$ is also Lagrangian.  This example signifies some 
  relationship of the content of this paper with the area of 
  integrable systems.

\
\noindent
{\bf (ii)} On the other hand if $X$ is a manifold equipped with an involute
  distribution $\cP \subseteq T_{X}$ of half the dimension of $X$ and $\cJ$
  is the generalized almost complex structure on $X$ corresponding to 
  an almost complex structure $J$ then $\cP$ and $\cJ$ are compatible if and 
  only if $J \cP \cap \cP = (0)$.  In other words the leaves of the foliation 
  are totally real sub-manifolds \cite{Bog}.  In this case the $\cJ$-compliment
  to $\cP$ is fixed uniquely as $J \cP$.     

\
\noindent
{\bf (iii)} One of the main classes of examples in this paper is where 
  $X$ is an $n$-torus bundle over an $n$-manifold and $\cJ$ is 
  a semi-flat generalized complex structure, $\cP$ is the vertical 
  foliation tangent to the torus fibers, and $\cQ$ is the horizontal 
  foliation given by the splitting of the tangent sequence given by 
  the connection as in section \ref{s:tbundle}. 

\end{examples} 

\begin{rem}
In the above and in much of what follows, the fact that $\cP$ is involute is
irrelevant.  That is to say, it could just be a sub-bundle of the 
tangent bundle.  However, it will be taken to be involute for the applications
that we have in mind, for instance when $\cP$ represents the tangent 
directions to a torus fibration.
\end{rem}

\begin{defin} \label{def:kahcompatible}
Suppose that $\cJ$ and $\cJ'$ constitute a generalized almost 
K\"{a}hler pair of 
generalized almost complex structures and $\cP \subseteq T_X$ is a sub-bundle 
of the tangent bundle of half the dimension.  Then we say that $\cP$ is 
compatible with the pair $(\cJ, \cJ')$ when 
\begin{eqnarray}
\cJ_2 (\Ann(\cP)) \subseteq \cP\\
\cJ_3 (\cP) \subseteq \Ann(\cP)\\
\cJ'_2 (\Ann(\cP)) \subseteq \cP\\
\cJ'_3 (\cP) \subseteq \Ann(\cP)\\
\cJ_1 \cJ'_1 (\cP) \subseteq \cP\\
\cJ'_1 \cJ_1 (\cP) \subseteq \cP\\
\cJ_4 \cJ'_4 (\Ann(\cP)) \subseteq \Ann(\cP)\\
\cJ'_4 \cJ_4 (\Ann(\cP)) \subseteq \Ann(\cP)
\end{eqnarray}
\end{defin}

Notice that if there is a sub-bundle $\cQ$ which is both a $\cJ$-compliment
and an $\cJ'$-compliment to $\cP$ then $\cP$ is compatible 
with $(\cJ, \cJ')$.  The converse will be shown below.  In the 
ordinary K\"{a}hler case \ref{def:genkahler} the 
condition that
$\cP$ is compatible with $(\cJ,\cJ')$ simply says that $\cP$ is Lagrangian 
with respect to the symplectic structure.
In the $B$-transformed almost K\"{a}hler case where we have an 
almost K\"{a}hler pair $(J, \om)$ 
the conditions are as follows: $\cP$
must be Lagrangian with respect to the symplectic structure $\om$ and also 
$B(\om^{-1} B \cP, \cP) = 0$ and  $B(J\cP, \cP) = 0$.

\begin{thm} \label{thm:compatabilities}
If $X$ is a $2n-$dimensional real manifold then a rank 
$n$ bundle $\cP \subseteq T_X$ is compatible with a 
generalized almost K\"{a}hler pair
$(\cJ, \cJ')$ if and only if there is a sub-bundle 
$\cQ \subseteq T_X$
 which is both a 
$\cJ$-compliment and a $\cJ'$-compliment to $\cP$.  These properties
specify $\cQ$ uniquely. 
\end{thm}

\noindent
{\bf Proof.}

If such a $\cQ$ exists then it is clear that $\cP$ and $\cQ$ 
are both compatible with the pair $(\cJ, \cJ')$.  Furthermore, the property
that $\cQ$ is a $\cJ$-compliment and a $\cJ'$-compliment to $\cP$ for
a generalized K\"{a}hler pair $(\cJ, \cJ')$ fixes $\cQ$ uniquely. Indeed, if 
we are in this situation and 
$G= -\cJ \cJ'$ is the generalized K\"{a}hler
metric then we have that 

\[
g - bg^{-1}b =G_3= -\cJ'_3 \cJ_1 - \cJ'_4 \cJ_3.
\]
  
Therefore the isomorphism $(g - bg^{-1}b)$ takes $\cP$ to 
$\Ann(\cQ)$ and so $\cQ$ must be the perpendicular to complement 
of $\cP$ with respect to the metric  $(g - bg^{-1}b)$.  We will now realize
$\cQ$ explicitly as the image of a different automorphism 
$K_+ \in GL(T_X)$ of the
tangent bundle with itself which becomes an almost complex structure 
in the when $b=0$.  Therefore the proof will be completed via the 
following lemma. 

\begin{lem} \label{lem:genkahler}

If $(\cJ, \cJ')$ is a generalized almost K\"{a}hler pair then
the map 
\[
K_{+} = J_{+} (1 - g^{-1} b) = \cJ_1 + \cJ'_1
\]
 is an isomorphism 
of the tangent bundle with itself.
If $\cP$ is compatible with the pair $(\cJ, \cJ')$ then $K_{+}$ 
takes $\cP$ to a sub-bundle $\cQ$, transversal to $\cP$,
and we have that  
 
\[
\cJ_{\big\vert \cP\oplus \op{Ann}(\cQ)} \; : \cP\oplus \op{Ann}(\cQ) \to \cQ
\oplus \Ann(\cP)
\]
\noindent
and
\[
\cJ'_{\big\vert \cP\oplus \op{Ann}(\cQ)} \; : \cP\oplus \op{Ann}(\cQ) \to \cQ
\oplus \Ann(\cP)
\]
\noindent
are isomorphisms of vector bundles.  In other words, $\cQ$ is both a 
$\cJ$-compliment and a $\cJ'$-compliment to $\cP$.

\end{lem}

\noindent
{\bf Proof.}  

\noindent
First of all, notice that $K_{+}$ is an isomorphism of the vector bundle
$T_X$ with itself.  Indeed, the vector bundle map $(g - b g^{-1} b)$ from
$T_X$ to $T_X^{\vee}$ corresponds to the metric $g(v,w) + g^{-1}(bv,bw)$ which
is positive definite and hence if we consider $f$ to $(g - b g^{-1} b)^{-1}$,
we see that 
\[
K_{+}(-(g+b)f J_{+})=1.  
\]

Now we have
\[
G_3 \cJ_1 =  -(\cJ_3 \cJ'_1 + \cJ_4 \cJ'_3)\cJ_1 -\cJ_3 \cJ'_1 \cJ_1 -\cJ_4 \cJ'_3 \cJ_1
\]

\[
=-\cJ_3 \cJ'_1 \cJ_1- \cJ_4(\cJ_3 \cJ'_1 + \cJ_4 \cJ'_3 - \cJ'_4 \cJ_3)
\]

\[
=-\cJ_3 \cJ'_1 \cJ_1 - \cJ_4 \cJ_3 \cJ'_1 - (\cJ_4)^2 \cJ'_3 + \cJ_4 \cJ'_4 \cJ_3
\]

\[
=-\cJ_3 \cJ'_1 \cJ_1 - \cJ_4 \cJ_3 \cJ'_1 - (-1 - \cJ_3 \cJ_2) \cJ'_3 + \cJ_4 \cJ'_4 \cJ_3
\]

\[
=-\cJ_3 \cJ'_1 \cJ_1 - \cJ_4 \cJ_3 \cJ'_1 +\cJ'_3  + \cJ_3 \cJ_2 \cJ'_3 + \cJ_4 \cJ'_4 \cJ_3
\]

\[
=-\cJ_3 \cJ'_1 \cJ_1 + \cJ_3 \cJ_1 \cJ'_1 +\cJ'_3  + \cJ_3 \cJ_2 \cJ'_3 + \cJ_4 \cJ'_4 \cJ_3.
\]

By inspection of the definition of compatibility of $\cP$ with the 
pair $(\cJ,\cJ')$
we have that all of these terms send $\cP$ into $\Ann(\cP)$.  Thus
$G_3 \cJ_1 (\cP) \subseteq \Ann(\cP)$.
Since the roles of $\cJ$ and $\cJ'$ are interchangeable we have,

\[
G_3 \cJ'_1 = 
-\cJ'_3 \cJ_1 \cJ'_1 + \cJ'_3 \cJ'_1 \cJ_1 +\cJ_3  + \cJ'_3 \cJ'_2 \cJ_3 + \cJ'_4 \cJ_4 \cJ'_3
\]
\noindent
as well.  Hence $G_3 \cJ'_1 (\cP) \subseteq \Ann(\cP)$.  
Therefore the image of 
$\cP$ under the isomorphism $K_+ = \cJ_1 + \cJ'_1$ is the perpendicular
sub-bundle to $\cP$ with respect to the metric $G_3$.

Now, in order to show the remaining 
claims, it suffices to define $\cQ$ as $K_{+}(\cP)$ and show that 
$\cJ_1(\cQ) \subseteq \cP$ and $\cJ_4 (\Ann(\cP)) \subseteq \Ann(\cQ)$.  
Indeed suppose that we have shown this.  Note that reversing the roles
of $\cJ$ and $\cJ'$ does not change $\cQ$ and so we get that 
$\cJ'_1(\cQ) \subseteq \cP$ and $\cJ'_4 (\Ann(\cP)) \subseteq \Ann(\cQ)$, and therefore
$\cJ(\cQ \oplus \Ann(\cP)) \subseteq \cP \oplus \Ann(\cQ)$ and 
$\cJ'(\cQ \oplus \Ann(\cP)) \subseteq \cP \oplus \Ann(\cQ)$, which is enough since
$\cJ$ and $\cJ'$ are isomorphisms.

Let $v$ be an element of a fiber of $\cQ$.  We may express it as 
$(\cJ_1 + \cJ'_1) w$ for a unique fiber $w$ of $\cP$ over the same point.  
Then 
\[
\cJ_1 v = {\cJ_1}^{2} w + \cJ_1 \cJ'_1 w 
= -w - \cJ_2 \cJ_3 w + \cJ_1 \cJ'_1 w
\]
which is an element of the fiber of $\cP$ over the same point.

Let $\mu$ be an element of a fiber of $\Ann(\cP)$.  Then, if $v$ is in the
fiber of $\cQ$ over the same point, we have that 
$(\cJ_4 \mu)v = -\mu (\cJ_1 v)$ which is 
zero by the previous paragraph.  Therefore $\cJ_4 \mu$ is in the 
fiber of $\Ann(\cQ)$ over the same point.     \hfill $\Box$   

\
 
The reader may wonder about the possibility of instead taking 
\[
K_{-} = \cJ_1 - \cJ'_1 = J_{-}(1 + g^{-1} b).
\] 

\begin{lem} \label{lem:allsame}

$K_{-}$ is an isomorphism of
the tangent bundle with itself.  In general it is not equal to $K_{+}$. 
However, if $\cP$ is compatible with the
generalized almost K\"{a}hler pair $(\cJ,\cJ')$, we have that 
$K_{+}(\cP) = K_{-}(\cP)$.  In fact we have that they are both equal to 
the orthogonal complement of $\cP$ with respect to the metric
$G_3 = g-bg^{-1}b$.
\end{lem}

\noindent
{\bf Proof.}
To see that $K_{-}$ is an isomorphism, simply note that 
$-f(g-b)J_{-}K_{-} = 1$ where $f$ is the inverse to $g-bg^{-1}b$.  
Define $\cQ_+ = K_{+}(\cP)$ and $\cQ_{-} = K_{-}(\cP)$. 
By the above arguments it is clear that $K_- = \cJ_1 - \cJ'_1$ is 
an isomorphism
from $\cP$ to the orthogonal complement of $\cP$ with respect to the metric 
$G_3 = g-bg^{-1}b$.  Therefore we have $K_{+}(\cP) = K_{-}(\cP)$. \hfill $\Box$

\begin{rem} \label{rem:partner}
As an aside, we mention that for any generalized almost 
complex structure, $\cJ$ there is another one $\cJ'$ such that
$(\cJ, \cJ')$ are a generalized almost K\"{a}hler structure.
\end{rem}

\noindent
Since we will not be using this and since the proof precisely mimics the proof
that every almost symplectic manifold has a compatible almost complex structure we do not include the proof here.

The next requirement that one should want to place on $(\cJ, \cJ', \cP)$ is 
that the distribution $\cQ$ be involute.  This is the analogue of considering 
a flat connection in definition \ref{def:semiflat}.  We plan to return to this 
analysis in a future paper.

\section{Examples}\label{s:examples}

\subsection{Mirror images of $B$-field and $\beta$-field
  transforms}\label{ss:Bmirror} 

\

Let $V$ be a vector bundle on $M$ with connection $\nabla$, $X = tot(V)$ and 
$\cJ=  F^{-1} (\pi^*\cK) F$ a generalized almost 
complex structure on $X$, where $F$ is defined in equation \ref{e:F}.  
We will need to relate 
$B$-field and $\beta$-field transforms of 
generalized complex structures $\cJ$ on $X$ to transformations of 
their mirror generalized complex structures 
$\cJwh = \Fwh^{-1} (\pih^*\cKwh) \Fwh$ on $\Xwh$.  

The transformation $\cK \to \exp(\underline{B}) \cK \exp(-\underline{B})$, where 

\begin{equation}\label{e:L}
\exp(\underline{B})=\left(
\begin{array}{cccc}

1    &    0    &    0   & 0 \\
0    &    1    &    0   & 0 \\
B_{31}  &   B_{32}   &    1   & 0 \\
B_{33}  &   B_{34}   &    0   & 1 

\end{array}
\right), \ \ \ \ \exp(\underline{B}) \in GL(V \oplus T_M \oplus V^{\vee} \oplus T_M^{\vee})
\end{equation}

\noindent
corresponds under mirror symmetry to the 
transformation $\cKwh \to \exp(\underline{\Bwh}) \cKwh \exp(-\underline{\Bwh})$, where

\begin{equation}\label{e:mirrorB}
\exp(\underline{\Bwh}) =\left(
\begin{array}{cccc}

1    &    B_{32}    &      B_{31}    & 0 \\
0    &    1         &      0         & 0 \\
0    &    0         &      1         & 0  \\
0    &    B_{34}    &      B_{33}    & 1 

\end{array}
\right), \ \ \ \ \exp(\underline{\Bwh}) \in GL(V^{\vee} \oplus T_M \oplus V \oplus T_M^{\vee}).
\end{equation}

\noindent
Similarly, the transformation $\cK \to 
\exp(\underline{\beta}) \cK \exp(-\underline{\beta})$, where

\begin{equation}\label{e:Q}
\exp(\underline{\beta})=\left(
\begin{array}{cccc}

1    &    0    &    \beta_{21}   & \beta_{22} \\
0    &    1    &    \beta_{23}   & \beta_{24} \\
0    &    0    &      1          & 0       \\
0    &    0    &      0          & 1 

\end{array}
\right), \ \ \ \ \exp(\underline{\beta}) \in  GL(V \oplus T_M \oplus V^{\vee} \oplus T_M^{\vee}).
\end{equation}

\noindent
corresponds under mirror symmetry to the 
transformation  $\cKwh \to \exp(\underline{\Betawh}) 
\cKwh \exp(-\underline{\Betawh})$, where

\begin{equation}\label{e:mirrorQ}
\exp(\underline{\Betawh})= \left(
\begin{array}{cccc}

1             &    0    &      0      & 0 \\
\beta_{23}    &    1    &      0      & \beta_{24} \\
\beta_{21}    &    0    &      1      & \beta_{22}       \\
0             &    0    &      0      & 1 

\end{array}
\right), \ \ \ \  \exp(\underline{\Betawh}) \in  GL(V^{\vee} \oplus T_M \oplus V \oplus T_M^{\vee}).
\end{equation}

\subsection{$B$-complex structures on $X=tot(T_M)$ and their mirror  
images}\label{ss:Bcplxmirror}  

\
Let us examine a very simple ``deformation'' of the setup from \cite{Pol}.  
It should be clear that there are many variants of this that one could easily 
do instead.  For instance one could vary the complex structure constructed on  
$TM$ from a fixed choice of connection.  Let 
$M$ be any manifold 
and $\nabla$ a flat and 
torsion-free connection 
on $T_M$ (one may drop the torsion free condtion, but then the analysis would become more complicated).  Let $X=tot(T_M)$, and $\Xwh =
tot(T_M^{\vee})$.  We will investigate $B$-field transforms of the canonical 
complex structure on $X$, where $B$ is an arbitrary real two-form.  We 
will give the condition for these transforms 
to be (integrable) $\nabla-$semi-flat (see definition \ref{def:semiflat})
 generalized complex structures on 
$X$ and give their integrable mirror structures on $\Xwh$.  

That is to say, consider a generalized almost complex
structure on $X$ of the form.  The $B-$field transform of the canonical 
complex structure is $\cJ=  F^{-1} (\pi^*\cK) F$ where 

\begin{equation}\label{e:K Bdom red}
\cK=\left(
\begin{array}{cccc}
0                 &         1            & 0   &   0  \\
-1                &         0            & 0   &   0  \\
-B_{32} - B_{33}        &      B_{31} - B_{34}       & 0   &   1   \\  
B_{31} - B_{34}         &      B_{32} + B_{33}       & -1  &   0  
\end{array}
\right).
\end{equation}

\noindent
and $B_{31}$ and $B_{34}$ represent arbitrary two forms on $M$, 
$B_{31} = -B_{31}^{\vee}$ 
and $B_{34} = -B_{34}^{\vee}$.  For this to be semi-flat, we 
need it to be adapted (to the splitting, see \ref{def:nablaadapted}) and 
hence $B_{32} +B_{33} = 0$.  Therefore we consider

\begin{equation}\label{e:K Bdom red2}
\cK=\left(
\begin{array}{cccc}
0                 &         1            & 0   &   0  \\
-1                &         0            & 0   &   0  \\
0                 &      B_{31} - B_{34}       & 0   &   1   \\  
B_{31} - B_{34}         &         0            & -1  &   0  
\end{array}
\right).
\end{equation}

Now the analysis in section \ref{s:integrability} tells us precisely when
the generalized almost complex structure $\cJ$ is 
integrable.  Namely, we must have that all Courant Brackets
of sections in the image of the subsheaf $\cS$ of flat sections of 
$T_M \oplus T_M^{\vee}$ under 

\begin{equation}\label{e:MBDom}
\cM=\left(
\begin{array}{cc}
-1              &    0 \\
B_{31} - B_{34}       &    1 
\end{array}
\right)
\end{equation}
\noindent
must vanish.  This, in turn, is equivalent to the following three Courant 
Brackets vanishing for any choice of flat sections $X, Y$  of $T_M$, 
and flat sections $\xi, \eta$  of $T_M^{\vee}$.
\[
\begin{split}
& \bigl[ -X + (B_{31} -B_{34})X,  -Y + (B_{31} -B_{34})Y \bigr]  = 0 \\
& \bigl[ -X + (B_{31} -B_{34})X, \eta \bigr]  = 0 \\
& \bigl[ \xi, \eta \bigr]  = 0
\end{split}
\]

Notice that when $B=0$, we recover no furthur conditions as expected.  
The second and third conditions are clearly vacuous.  
Set $B' =B_{31} - B_{34}$.  The first condition then reads:
\[
\myiota_{-X} d(B'Y) - \myiota_{-Y}d(B'X) 
- \frac{1}{2}d(\myiota_{-Y}(B'X) - \myiota_{-X}(B'Y)) = 0
\]
or 
\[
\begin{split}
0 & = -\myiota_{X}d\myiota_{Y}B' + -\myiota_{Y}d\myiota_{X}B' +
d(\myiota_{Y}\myiota_{X}B') \\
& = -\myiota_{X}(\cL_{Y} -
\myiota_{Y}d)B'+\myiota_{Y}(\cL_{X} - \myiota_{X}d)B' + (\cL_{Y} -
\myiota_{Y}d)\myiota_{X}B' \\
& = -\myiota_{[X,Y]}B' +
2\myiota_{X}\myiota_{Y}dB' +\cL_Y \myiota_X B'-\myiota_Y \cL_X B' -
\myiota_Y \myiota_X dB'\\
&  = 3\myiota_X \myiota_Y dB'. 
\end{split}
\]

\noindent
Thus $\cJ$ is integrable if and only if $B' = B_{31} - B_{34}$ is closed.  
The mirror structure on $\Xwh$ is given by  
$\cJwh = \Fiwh (\pih^*{\cKwh}) \Fwh$
where 

\begin{equation}\label{e:Dom cKwhB}
{\cKwh}=\left(
\begin{array}{cccc}

0    & B_{31}-B_{34}    & 0              & 1   \\
0    & 0          & -1             & 0     \\
0    & 1          & 0              & 0     \\
-1   & 0          & B_{31}-B_{34}        & 0

\end{array}
\right)
\end{equation}

\noindent
This is the $\beta$-field transform of the canonical symplectic structure 
on $\Xwh$, where 

\[ \beta = 
\left(\begin{array}{cc} \jh & \alh \end{array}\right)
\matr{0}{\pih^*(B_{31} - B_{34})}{0}{0}
\left(\begin{array}{c} {\jh}^{\vee} \\ {\alh}^{\vee} \end{array}\right)
= \jh \pih^*(B_{31} - B_{34}) {\alh}^{\vee}
\]

\sbr

\subsection{$B$-symplectic structures on $\Xwh = tot({T_M}^{\vee})$ and 
their mirror transforms}\label{ss:Bsympmirror} 

\

Let $\nabla^{\vee}$ be the dual of a flat, torsion-free connection 
$\nabla$ on $T_M$.  In this section we will compute the conditions for
a $B$-field transform of the canonical symplectic structure on $\Xwh$ to be 
a $\nabla^{\vee}$-semi-flat generalized complex structure and 
find the (integrable) mirror image structure on $X$.

This $B$-symplectic generalized almost complex structure  
$\cJwh = \Fwh^{-1} \cKwh \Fwh$ on $\Xwh$ is given by 

\begin{equation}\label{e:K Bcan red}
{\cKwh}=\left(
\begin{array}{cccc}
-B_{33}                    &         -B_{34}            & 0      &   1  \\
B_{31}                     &         B_{32}             & -1     &   0  \\
B_{32}B_{31}+B_{31}B_{33}  &  1+(B_{32})^2 -B_{31}B_{34}& -B_{33}    B_{31}\\  
-1+B_{34}B_{31}-(B_{33})^2 & B_{34}B_{32}- B_{33}B_{34} &  -B_{34}  &  B_{33}  
\end{array}
\right).
\end{equation}
\noindent
The adapted requirement (see \ref{def:nablaadapted}) 
forces $B_{32} = B_{33} = 0$ and so

\begin{equation}\label{e:KBcanred2}
{\cKwh}=\left(
\begin{array}{cccc}
0            & -B_{34}            & 0      &   1  \\
B_{31}          & 0               & -1     &   0  \\
0            & 1- B_{31} B_{34}      & 0      &   B_{31}   \\  
-1+ B_{34} B_{31}  & 0               & -B_{34}   &   0  
\end{array}
\right).
\end{equation}

Now using the above analysis on integrability we know that this the 
generalized almost complex structure on $\Xwh$ will be 
integrable if and only if all Courant Brackets
of sections in the image of the subsheaf $\cS$ of flat sections of 
$T_M \oplus T_M^{\vee}$ under 

\begin{equation}\label{e:M B Dom}
\cM=\left(
\begin{array}{cc}
-1         &    -B_{31} \\
-B_{34}        &    1 - B_{34} B_{31} 
\end{array}
\right)
\end{equation}

\noindent
must vanish.  This, in turn, is equivalent to the following three Courant 
Brackets vanishing for any choice of $X, Y$ flat sections of $T_M$, 
and $\xi, \eta$ flat sections of $T_M^{\vee}$.

\[
\begin{split}
& \bigl[ -X  -B_{34} X,  -Y + -B_{34} Y \bigr]  = 0 \\
& \bigl[ -X  -B_{34} X, -B_{31} \eta + \eta - B_{34} B_{31} \eta \bigr]  = 0 \\
& \bigl[ -B_{31} \xi + \xi - B_{34} B_{31} \xi,  -B_{31} \eta + \eta - B_{34} B_{31} \eta \bigr]  = 0
\end{split}
\]

As in the previous subsection, the first equation is equivalent 
to $dB_{34} = 0$.  The second equation is equivalent to 

\[
\begin{split}
& [ X  ,  B_{31} \eta]  = 0 \\
& \myiota_X d \myiota_{B_{31} \eta} B_{34} - \myiota_{B_{31} \eta} d \myiota_{X} B_{34}
-\frac{1}{2} d( \myiota_{B_{31} \eta} \myiota_X B_{34} - \myiota_X (B_{34} B_{31} \eta)) = 0
\\
\end{split}
\]

The first of these equations simply says that $B_{31}$ is a flat bivector field.
On the other hand, we claim that if $dB_{34} = 0$ and $B_{31}$ is a flat 
bivector field then the second part of the second equation and also the 
third equation are also satisfied, and hence all the equations are satisfied.
Indeed, we have 

\[
d\eta = d\xi = d(\myiota_{B_{31} \xi} \eta) =  d(\myiota_{B_{31} \eta} \xi) = 0
\]

\noindent
Therefore the third equation gives

\[
\begin{split}
&\myiota_{B_{31} \xi} d \myiota_{B_{31} \eta} B_{34} 
-\myiota_{B_{31} \eta} d \myiota_{B_{31} \xi} B_{34}
+\frac{1}{2}(d \myiota_{B_{31} \xi}  \myiota_{B_{31} \eta} B_{34} 
-d \myiota_{B_{31} \eta} \myiota_{B_{31} \xi} B_{34}) \\
&=\myiota_{B_{31} \xi}\cL_{B_{31} \eta} B_{34} - \myiota_{B_{31} \eta}\cL_{B_{31} \xi} B_{34}
+d \myiota_{B_{31} \xi}  \myiota_{B_{31} \eta} B_{34} \\
&=\myiota_{B_{31} \xi}\cL_{B_{31} \eta} B_{34} - \myiota_{B_{31} \eta}\cL_{B_{31} \xi} B_{34}
-\myiota_{B_{31} \xi} d \myiota_{B_{31} \eta} B_{34} 
+\cL_{B_{31} \xi} \myiota_{B_{31} \eta} B_{34} \\
&=\myiota_{[B_{31} \xi, B_{31} \eta]} B_{34} +\myiota_{B_{31} \xi}\cL_{B_{31} \eta} B_{34}           - \myiota_{B_{31} \xi} d \myiota_{B_{31} \eta} B_{34}\\
&=\myiota_{B_{31} \xi}\cL_{B_{31} \eta} B_{34}
- \myiota_{B_{31} \xi} d \myiota_{B_{31} \eta} B_{34}
= \myiota_{B_{31} \xi}\cL_{B_{31} \eta} B_{34} 
- \myiota_{B_{31} \xi} \cL_{B_{31} \eta} B_{34} \\
&=0
\end{split}
\]

\noindent
Similarly the second part of the second equation gives 

\[
\begin{split}
& \myiota_X d \myiota_{B_{31} \eta} B_{34} - \myiota_{B_{31} \eta} d \myiota_{X} B_{34}
-\frac{1}{2} d( \myiota_{B_{31} \eta} \myiota_X B_{34} 
- \myiota_X \myiota_{B_{31} \eta}B_{34} )\\
& = \myiota_X \cL_{B_{31} \eta} B_{34} - \myiota_{B_{31} \eta} \cL_X B_{34}
-\cL_{B_{31} \eta} \myiota_X B_{34} + \myiota_{B_{31} \eta} d \myiota_X B_{34}\\
&=\myiota_{[X, B_{31} \eta]} B_{34} - \myiota_{B_{31} \eta} \myiota_X d B_{34}\\
&=0
\end{split}
\]

\noindent
The mirror structure ${\cK}$ to ${\cKwh}$ is thus given by

\begin{equation}\label{e:K Bcan red2}
{\cK}=\left(
\begin{array}{cccc}
0            & 1 -B_{31}B_{34}   & 0                 &   B_{31}  \\
-1           & 0                 & B_{31}            &   0  \\
0            & -B_{34}           & 0                 &   1   \\  
-B_{34}      & 0                 &-1 + B_{34}B_{31}  &   0  
\end{array}
\right).
\end{equation}

\noindent
Therefore the mirror $\cJ$ of $\cJwh$ is the canonical complex structure 
transformed by the composition of the $B$-field 
\[
\begin{split}
(d\pi)^{\vee} (\pi^*B_{34}) (d\pi) 
\end{split}
\]
\noindent
and the $\beta$-field  
\[
\begin{split}
j (\pi^*B_{31}) j^{\vee}.
\end{split}
\]

\medskip

\noindent
{\bf Acknowledgements}\label{ss:thanks} I would like to thank 
Tony Pantev for his guidance and encouragement and also Mitya Boyarchenko 
for previous collaboration.  Thanks to Alessandro Tomasiello  
and Sukhendu Mehrotra for helpful discussions.

\end{document}